\voffset=-0.6in
\magnification=\magstep1
\input amstex
\documentstyle{amsppt}
\magnification 1200
\catcode`\@=11
\def\logo@{\null}
\catcode`\@=12
\NoBlackBoxes

\def\ve{\varepsilon}

\def\R{\Bbb R}
\def\C{\Bbb C}

\def\P{\Cal P}

\def\Calg{$C^*$-algebras }
\def\Cag{$C^*$-algebra }
\def\Csalg{$C^*$-subalgebra }
\def\Csalgs{$C^*$-subalgebras }
\def\AW{$AW^*$-algebra }
\def\AWs{$AW^*$-algebras }
\def\sAW{$AW^*$-subalgebra }

\def\Wag{$W^*$-algebra }

\def\ve{\varepsilon}
\def\Re{\hbox{Re }}
\def\Im{\hbox{Im }}
\def\card{\hbox{card }}
\vglue 2 true cm
\centerline{\bf ABELIAN STRICT APPROXIMATION IN AW*-ALGEBRAS}
\centerline{\bf AND WEYL-VON NEUMANN TYPE THEOREMS \footnote""
}
\leftheadtext{C. D'Antoni and L. Zsid\`o}
\rightheadtext{Weyl-von Neumann type theorems}
\vglue 1 true cm
\centerline{Claudio D'Antoni and L\'aszl\'o Zsid\'o\footnote""{Supported
by  MIUR, INDAM and EU.}}
\vglue 1 true cm
\centerline{\it Dedicated to Professor E. Effros on his
$\, 70^{\text{th}}$ birthday}
\vglue 1.5 true cm

\topmatter
\abstract\nofrills
{\smc Abstract.} In this paper, for a \Cag $A$ with $M=M(A)$ an
$AW^*$-algebra, or equivalently, for an essential, norm-closed,
two-sided ideal $A$ of an \AW $M\,$, we investigate the strict
approximability of the elements of $M$ from commutative \Csalgs of
$A\,$. In the relevant case of the norm-closed linear span $A$
of all finite projections in a semi-finite \AW $M$ we shall give a
complete description of the strict closure in $M$ of any maximal
abelian self-adjoint subalgebra (masa) of $A\,$. We shall see that
the situation is completely different for discrete respectively
continuous $M\,$:

In the discrete case, for any masa $C$ of $A\,$, the strict
closure of $C$ is equal to the relative commutant $C'\cap M\,$,
while in the continuous case, under certain conditions concerning
the center valued quasitrace of the finite reduced algebras of $M$
(satisfied by all von Neumann algebras), $C$ is already strictly
closed. Thus in the continuous case no elements of $M$ which are
not already belonging to $A$ can be strictly approximated from
commutative $C^*$-subalgebras of $A\,$.

In spite of this pathology of the strict topology in the case of
the norm-closed linear span of all finite projections of a
continuous semi-finite $AW^*$-algebra, we shall prove that in
general situations including also this case, any normal $y\in M$
is equal modulo $A$ to some $x\in M$ which belongs to an
order theoretical closure of an appropriate commutative
$C^*$-subalgebra of $A\,$. In other words, if we replace the
strict topology with order theoretical approximation,
Weyl-von Neumann-Berg-Sikonia type theorems will hold in
substantially greater generality.
\endabstract
\endtopmatter

\smallskip

{\bf Introduction }
\bigskip

Let $A$ be a $C^*$-algebra. The {\it multiplier algebra} of $A$ is
the \Csalg
\medskip

\centerline{$\{x\in A^{**};\  xa,\  ax\in A\hbox{ for all } a\in A\}$}
\medskip

\noindent
of the second dual $A^{**}$ (see [Ped 2], Section 3.12 or [WO],
Chapter 2). A natural locally convex vector space topology on $M(A)\,$,
called the {\it strict topology} $\beta\,$, is defined by the seminorms
\medskip

\centerline{$x\mapsto\|xa\|\hbox{ and } x\mapsto\|ax\|\, ,\qquad a
\in A\, .$}
\medskip

\noindent
It is complete and compatible with the duality between $M(A)$ and $A^*$.
Hence the strict topology is weaker than the norm-topology on $M(A)$, 
but stronger than the restriction to $M(A)$ of the weak * topology of 
$A^{**}$.

We notice that for $A$ the \Cag $K(H)$ of all compact linear
operators on a complex Hilbert space $H\,$, $M(A)$ can be identified
with the \Cag $B(H)$ of all bounded linear operators on $H$ and
on every bounded subset of $B(H)$ the strict topology coincides with
the $s^*$-topology.

More generally, if $M$ is an \AW (see [Kap 1] or [Be], \S 4 or
[S-Z], \S 9) and $A$ is an essential, norm-closed, two-sided ideal
of $M\,$, then, by a theorem of B. E. Johnson, $M$ can be identified
with $M(A)$ (see [J] or [Ped 3]). Thus the pairs $\big( A\, ,\,
M(A)\big)\,$, where $A$ is a \Cag such that $M(A)$ is an
$AW^*$-algebra, are exactly the pairs $\big( A\, ,\, M\big)\,$,
where $M$ is an \AW and $A$ is an essential, norm-closed, two-sided
ideal of $M\,$.

A relevant case of essential, norm-closed, two-sided ideal of an
\AW is the norm-closed linear subspace $A$ generated by all finite
projections of a semi-finite $AW^*$-algebra $M\,$. Then there are
central projections $p_1\, ,\, p_2\, ,\, p_3$ of $M$ with
$p_1+p_2+p_3=1_M$ such that $M p_1$ is finite, $M p_2$ is properly
infinite and discrete, while $M p_3$ is properly infinite and
continuous (see [Be], \S 15, Theorem 1). Since $A p_1=M p_1\,$,
the non-trivial cases are $A p_2$ and $A p_3\,$, with $M(A p_2)
=M p_2$ properly infinite and discrete and $M(A p_3)=M p_3$ properly
infinite and continuous.

In the previous paper [D-Z] we investigated the strict approximability
of a normal element $x$ of $M(A)$ from a commutative \Csalg of $A\,$.
More precisely, we say that $x$ belongs to the {\it abelian strict
closure} of $A$ if there exists a commutative \Csalg $C_x$ of $A$
such that $x\in\overline{C_x}^{\,\beta}$. Abelian strict
approximability is closely related to the classical
Weyl-von Neumann-Berg-Sikonia (WNBS) Theorem, which claims that
in the case of $A=K(H)\,$, $H$ a separable complex Hilbert space,
every normal element of $M(A)=B(H)$ is of the form $a+x$ with
$a\in A$ and $x$ in the abelian strict closure of $A\,$.
For a general $\sigma$-unital \Cag $A\,$, that is a \Cag having
a countable approximate unit, we proved a partial extension in
[D-Z], Theorem 1, which implies that all elements $y\in M(A)$ are
of the form $a+x_1+x_2\,$, where $a\in A$ and $x_1\in B_1\, ,\,
x_2\in B_2\,$, where $B_1\, ,\, B_2$ are separable \Csalgs of
$M(A)$ such that every normal element of $B_j\, ,\, j=1,2\,$,
belongs to the abelian strict closure of $A\,$. Moreover, if $y$
is self-adjoint then $x_1\, ,\, x_2$ can be chosen self-adjoint,
so in this situation $x_1\, ,\, x_2$ themselves belong to the
abelian strict closure of $A\,$.

We notice that if the multiplier algebra of a $\sigma$-unital
$C^*$-algebra $A$ is of real rank zero (see [Br-Ped]), then,
according to [M] and [Zh], the WNBS Theorem holds in the same
formulation as in the classical case.

In this paper we discuss abelian strict approximability for
a \Cag $A$ which is the norm-closed linear subspace generated
by all finite projections of some semi-finite \AW $M\,$.
Since the abelian strict closure of $A$ is the union of all
$\overline{C}^{\,\beta}$ with $C$ a maximal abelian self-adjoint
subalgebra ({\it masa}) of $A\,$, we are interested in describing
$\overline{C}^{\,\beta}$ for any masa $C$ of $A\,$. We shall
see that the situation is completely different for discrete
respectively continuous $M(A)=M\,$:

In the discrete case $\,\overline{C}^{\,\beta}$ is equal to the
relative commutant $C'\cap M(A)$ (Theorem 1), while in the
continuous case, under a certain condition on the centre
valued quasitrace of the reduced \sAW of $M(A)$ by a finite
projection of central support $1_{M(A)}$ (always satisfied if
$M(A)$ is a von Neumann algebra), $C$ is already strictly
closed (Theorem 3).

Consequently, if $M(A)$ is a properly infinite, continuous \AW
satisfying the above mentioned condition, then the unit of $M(A)$
does not belong to the abelian strict closure of $A\,$, that is
it does not exist an approximate unit for $A$ contained in a
commutative $*$-subalgebra of $A\,$. In particular, in this case
$A$ is not $\sigma$-unital. We notice that it was already shown
in [Ak-Ped], Proposition 4.5, that the norm-closed linear span of
all finite projections of a type II${}_\infty$ factor is a
non-$\sigma$-unital $C^*$-algebra. Nevertheless, also in this
case WNBS type theorems can be proved.
Indeed, if $A$ is the norm-closed linear subspace generated by
all finite projections of some countably decomposable semi-finite
$W^*$-factor $M\,$, then, according to [Z], Theorem 3.1, every
normal $y\in M(A)=M$ is of the form $a+x$ with $a\in A$ and $x$
in the $s^*$-closure in $M$ of some masa $C$ of $A\,$. Since
the $s^*$-closure of a commutative $*$-subalgebra of a
$W^*$-algebra is equal to its monotone order closure (cf.
[Kad 1] and [Ped 1]), it is natural to expect that for extensions
of the WNBS Theorem to non-$\sigma$-unital $C^*$-algebras the
strict closure should be replaced by an order theoretical closure.
Along this line we prove several WNBS type theorems in a general
setting which includes the case of the norm-closed linear span
of all finite projections of a countably decomposable semi-finite
$AW^*$-algebra.

More precisely, we prove that if $\Cal J$ is a norm-closed
two-sided ideal of a (unital) Rickart $C^*$-algebra $M$, which
has a countable ``order theoretical approximate unit'', then
any normal $y\in M$ is of the form $y=a+x\,$, where $a\in A$
is of arbitrarily small norm and $x$ belongs to the order
theoretical closure of some masa of $\Cal J$ (Theorem 4 and the
subsequent remark). Moreover, the above $x$ can be chosen as a
particular infinite linear combination of a sequence of mutually
orthogonal projections from $\Cal J$ (Theorems 5 and 6).

Since only little of the specific properties of Rickart
$C^*$-algebras is used, we are left with the question, to which
extent the above mentioned WNBS type theorems hold if $M$ is
assumed to be only a $C^*$-algebra of real rank zero.
\bigskip


{\bf 1 Abelian Strict Closure in Discrete AW*-algebras }
\bigskip

First we prove a general result concerning a masa $\, C$ of a
$C^*$-algebra $A\,$, whose multiplier algebra is an
$AW^*$-algebra, that is, according to the theorem of B. E. Johnson
quoted in Introduction (see [J] or [Ped 3]), a masa $\, C$ of an
essential, norm-closed, two-sided ideal $A$ of some $AW^*$-algebra.
We notice that a part of this result holds for a masa of an
essential, norm-closed, two-sided ideal of any Rickart $C^*$-algebra.
We shall restrict us to unital Rickart $C^*$-algebras, because
adjoining a unit to a non-unital Rickart $C^*$-algebra $M$, we
obtain a unital Rickart $C^*$-algebra $\widetilde M$ (see
[Be], \S 5, Theorem 1 or [S-Z], 9.11.(1)) and it is easy to see
that every essential, norm-closed, two-sided ideal of $M$ is an
essential, norm-closed, two-sided ideal also of $\widetilde M\,$.

Any essential two-sided ideal $\Cal J$ of a \Cag $M$ induces
a strict topology $\beta_{\Cal J}$ on $M\,$, defined by the
seminorms
\smallskip

\centerline{$M\ni x\mapsto\|xa\|\hbox{ and } x\mapsto\|ax\|\, ,
\qquad a\in \Cal J\, .$}
\smallskip

\noindent With this definition, the usual strict topology 
on the multiplier algebra of a \Cag $A$ is $\beta_A\,$.

For the basic facts concerning Rickart \Calg and \AWs see [Be],
\S\S$\,$ 3, 4 and 5, or [S-Z], \S 9.
\medskip

\proclaim{Lemma 1} Let $M$ be a unital $\, C^*$-algebra, $\Cal J$ an
essential, norm-closed, two-sided ideal of $\, M$, and $\, C$ a
masa of $\Cal J\,$. By the strict topology on $M$ we shall
understand $\beta_{\Cal J}\,$, which of course is the usual strict
topology when $M$ is an $AW^*$-algebra and so can be identified
with the multiplier algebra $M(\Cal J)\,$. Then

\itemitem{{\rm (i)}} every $x\ge 0$ in the strict closure of
$\, C$ in $M$ belongs to the strict closure of
$\{ b\in C\, ; 0\le b\le x\}$ in $M\,$.
\smallskip

\noindent Let us next assume that $M$ is a Rickart $\, C^*$-algebra.
Then
\smallskip

\itemitem{{\rm (ii)}} for every $0\le b\in C$ and every $\delta >0$
there is a  projection $f_\delta\in C$ such that
\medskip

\centerline{$bf_\delta\ge \delta f_\delta,\quad b(1_M -f_\delta)
\leq\delta(1_M-f_\delta)\, ,$} 
\medskip

\quad\, so $\, C$ is the norm-closed linear span of its projections$;$

\itemitem{{\rm (iii)}} any projection $e$ in the strict closure of
$\, C$ in $M$  belongs to the strict closure of $\{f\in C\, ; f\le
e\hbox{ projection }\}$ in $M\, ;$

\itemitem{{\rm (iv)}} any projection $e$ in the relative commutant
$\, C'\cap M$ is the least upper bound of $\{f\in C\, ; f\le
e\hbox{ projection }\}$ in the projection lattice of $M\,$,
in particular $\, C'\cap M$ is a masa of $M\,$.
\smallskip

\noindent Finally, assuming $M$ to be an $AW^*$-algebra,
\smallskip

\itemitem{{\rm (v)}} the relative commutant $\, C'\cap M$ is
the \sAW of $M$ generated by $\, C\,$, so $\, C'\cap M$ can be
identified with $M(C)\, ;$

\itemitem{{\rm (vi)}} the strict closure of $\, C$ in $M$
coincides with $\, C'\cap M$ if and only if $\, C$ contains a
two-sided approximate unit  for $\Cal J\,$, in which case the
strict topology of $\, M(C)=C'\cap M$ is the restriction of the
strict topology of $\, M(\Cal J)=M\,$.
\endproclaim

{\bf Proof.} The strict closure $\,{\overline{C}}^{\beta_{\Cal J}}$
of $\, C$ being an abelian \Csalg of $M(A)\,$, we have for every
$b\in C$
\medskip
\centerline{$(x-b)^*(x-b)\ge (x-\Re b)^2\ge\big( x-(\Re b)_+
\big)^2\ge (x-b_o)^2\, ,$}
\medskip

\noindent where
\smallskip

\centerline{$b_o=\frac12\Big( x+(\Re b)_+-|x-(\Re b)_+|\Big)$}
\medskip

\noindent
denotes the greatest lower bound of $x$ and $(\Re b)_+$ in the
Hermitian part of ${\overline{C}}^{\beta_{\Cal J}}\, .$ Since 
\smallskip

\centerline{$0\le b_o\le (\Re b)_+\in C\subset \Cal J\, ,$}
\medskip

\noindent
by [Ped 2], Prop. 1.4.5 we have $b_o\in \Cal J\, ,$ so
\medskip

\centerline{$b_o\in C'\cap \Cal J=C\, .$}
\medskip

\noindent
Thus, for every $a\in \Cal J$ and $b\in C$ we have $\| (x-b)a\|\ge
\|(x-b_o)a\|$ for some $0\le b_o\le x$ in $C$ and (i) follows.

For (ii) put
\smallskip

\centerline{$f_\delta=\hbox{ support of } (b-\delta 1_{A^{**}})_+
\hbox{ in } M\, .$}
\medskip

\noindent
Then $f_\delta$ commutes with every element of $C$ and
\medskip

\centerline{$bf_\delta\ge \delta f_\delta,\quad b(1_M-f_\delta)
\le\delta(1_M -f_\delta)\, .$}
\medskip

\noindent
In particular, $f_\delta\le\frac1\delta b\in A$ and [Ped 2],
Prop.1.4.5 yields $f_\delta\in \Cal J\, .$ Consequently
$f_\delta\in C'\cap A=C\, .$

For (iii) let $0\ne a\in \Cal J$ and $\ve > 0$ be arbitrary.
According to (i) there exists $0\le b\le e$ in $C$ such that 
\smallskip

\centerline{$\displaystyle \|(e-b)a\| <\frac\ve2\, .$}
\medskip

\noindent
Further, by (ii) there is a projection $f\in C$ with
\medskip

\centerline{$\displaystyle bf\ge\frac\ve{2\|a\|}f,\qquad
b(1_{A^{**}}-f)\le\frac\ve{2\|a\|}\cdot (1_{A^{**}}-f)\, .$}
\medskip

\noindent
Then $f\le e$ and $e-f\le(e-bf)^2,$ so
$$\aligned
\|(e-f)a\|&=\|a^*(e-f)a\|^{1/2}\le\\
&\le\|a^*(e-bf)^2a\|^{1/2}=\|(e-bf)e\|\le\\
&\le\|(e-b)e\|+\| b(1_{A^{**}}-f)e\|<\\
&<\frac\ve2+\frac\ve{2\|a\|}\|a\|=\ve
\endaligned
$$

For (iv) we have to show that if a projection $g\in M$
majorizes all projections $C\ni f\le e\,$, then $g\ge e\,$,
that is $e$ is equal to the greatest lower bound $e\land g$
of $e$ and $g$ in the projection lattice of $M\, .$
Let us assume that
\medskip

\centerline{$e_o=e-e\land g\ne 0\, .$}
\medskip

\noindent
Since $\Cal J$ is essential ideal in $M\,$, there exists
$a\in \Cal J$ with $ae_o\ne 0\,$. Choosing some $0<\delta<
\|e_oa^*ae_o\|$ and putting
\medskip

\centerline{$e_1=$ support of $(e_oa^*ae_o-\delta 1_M)_{+}$
in $M\,$,}
\smallskip

\noindent we have 
\smallskip

\centerline{$0\ne e_1 \le\frac1\delta e_oa^* ae_o\in \Cal J\, .$}
\medskip

\noindent
Clearly, $e_1\le e_o$ and [Ped 2], Prop. 1.4.5 yields also
$e_1\in \Cal J\, .$ Furthermore, for every projection $f\in C$
we get successively
$$
\aligned
f e\,&\in C'\cap \Cal J=C\hbox{ and } fe\le e,\\
f e\,&\le e\land g, \hbox{ hence } fe_o=(fe)e_o=0,\\
fe_1&= (fe_o)e_1=0\, .
\endaligned
$$
\noindent
Taking into account (ii), it follows that
\medskip

\centerline{$be_1=0\,\hbox{ for all } b\in C\, ,$}
\medskip

\noindent in particular
\medskip

\centerline{$e_1\in C'\cap \Cal J=C\, .$}
\medskip

\noindent But then $e_1\le e_o\le e$ implies $e_1\le e\land g\,$,
which contradicts $0\ne e_1\le e_o=e-e\land g\, .$

In particular, $C'\cap M$ is commutative. For the proof we notice
that, since $C'\cap M$ is a Rickart $C^*$-subalgebra of $M$ (see
[Be], \S 5, Proposition 5 or [S-Z], 9.12.(1)), it is the
norm-closed linear span of its projections (see e.g. [S-Z], 9.4)
and therefore it is enough to show that any two projections
$e_1 , e_2\in C'\cap M$ commute. But the $*$-automorphism
$M\ni x\longmapsto (2\;\! e_2-1_M) x(2\;\! e_2-1_M)\in M$ leaves
fixed $C\,$, hence also the least upper bound of any projection
family in $C$ in the projection lattice of $M\,$. Therefore it
leaves fixed $e_1\,$, that is $e_1 e_2 =e_2 e_1\,$.

Moreover, $C'\cap M$ is a masa of $M\,$. Indeed, if $C_o\supset
C'\cap M$ is a commutative subalgebra of $M\,$, then $C_o\supset C$
and thus we have also $C_o\subset C_o{}'\cap M\subset C'\cap M\,$.

For (v) we first notice that $C'\cap M$ is an \sAW of $M$ containing
$C$ (see [Be], \S 4, Proposition 8 or [S-Z], 9.24.(1)). Now let $N$
be any \sAW of $M$ containing $C\,$. By (iv) $N$ contains all
projections from $C'\cap M\,$, hence $N\supset C'\cap M\,$.
Consequently $C'\cap M$ is the \sAW of $M$ generated by $C\,$.

Further, $C$ is a two-sided ideal of $C'\cap M:$
\medskip

\centerline{$b\in C\hbox{ and } y\in C'\cap M\;\Longrightarrow
\; b y\in C'\cap \Cal J=C\, .$}
\medskip

\noindent
Moreover, it is essential, because a projection $e\in C'\cap M$ with
$C e=\{0\}$ belongs to the \sAW of $C'\cap M$ generated by $C$ only
if $e=0\,$. Hence we can identify $C'\cap M$ with $M(C)$ (see [J] or 
[Ped 3]).

Finally we prove (vi). If the strict closure of $C$ in $M$ is
$C'\cap M\ni 1_M\, ,$ then there exists a net $(u_\iota)_\iota $ in
$C$ with $u_\iota \to 1_M$ strictly in $M\,$, that is
\medskip

\centerline{$\|a-u_\iota  a\|\to 0\hbox{ and } \|a-u_\iota  a\|
\to 0\hbox{ for all }a\in \Cal J\, .$}
\medskip

Conversely, let us assume that $C$ contains a two-sided approximate
unit $(u_\iota )_\iota $ for $\Cal J\, .$ Then the strict topology
$\beta_C$ of $M(C)=C'\cap M$ agrees with the strict topology
$\beta_{\Cal J}$ of $M(\Cal J)=M$ on every norm bounded subset of
$C'\cap M\, .$ Indeed, if  $(y_\lambda)_\lambda$ is a norm bounded
net in $C'\cap M\,$, convergent to 0 with respect to $\beta_C\,$,
and $0\ne a\in \Cal J ,\, \ve> 0$ are arbitrary, then there exists
$\iota_o$ such that
\medskip
\centerline{$\displaystyle \|y_\lambda\|\cdot\|a-u_{\iota _o}a\|
<\frac\ve2\hbox{ for all } \lambda\, ,$}
\medskip

\noindent and then there exists some $\lambda_o$ with
\medskip
\centerline{$\displaystyle \|y_\lambda u_{\iota _o}\|<
\frac{\ve}{2\|a\|}\hbox{ for every }\lambda\ge\lambda_o\, .$}
\medskip

\noindent It follows for every $\lambda\ge\lambda_o :$
\medskip
\centerline{$\displaystyle \|y_\lambda a\|\le\|y_\lambda
(a-u_{\iota _o}a)\|+\|y_\lambda u_{\iota_o}a\| <\frac\ve2+
\frac\ve{2\|a\|}\|a\|=\ve\, .$}
\medskip

\noindent But $\beta_C$ is the finest locally convex vector space
topology on $C'\cap M$ that agrees with $\beta_C$ on every norm
bounded subset of $C'\cap M$ (see [T], Cor. 2.7). Thus the
restriction of $\beta_{\Cal J}$ to $C'\cap M\,$, which is plainly
finer than $\beta_C\,$, is actually equal to $\beta_C\,$. In
particular, the $\beta_C$-density of $C$ in $M(C)$ implies the
$\beta_{\Cal J}$-density of $C$ in $C'\cap M\,$.

\hfill$\square$
\bigskip

It is well known that every commutative \AW $Z$ is monotone complete
(see e.g. [S-Z], 9.26, Proposition 1). If $M$ is an arbitrary
$AW^*$-algebra, we call
\medskip
\centerline{$\Phi:\big\{ e\in M; e\hbox{ projection }\big\} \to Z^+$}
\medskip

\noindent completely additive whenever, for every family
$(e_\iota )_\iota $ of mutually orthogonal projections in $M,$ we have
\smallskip
\centerline{$\displaystyle \Phi\big(\bigvee_\iota e_\iota\big) =
\sum_\iota \Phi(e_\iota )\, ,$}
\smallskip

\noindent where the sum stands for the least upper bound in $Z^+$
of all finite sums of $\Phi(e_\iota )$'s.
\smallskip

Now we describe the strict closure of a masa of the norm-closed
two-sided ideal generated by the finite projections of a discrete
semi-finite \AW :
\medskip

\proclaim{Theorem 1 {\rm $($on the abelian strict closure in discrete
$AW^*$-algebras$)$}}
Let $M$ be a discrete $AW^*$-algebra, $A$ the norm-closed
linear span of all finite projections of $M\,$, and $\, C$ a masa
of $A\,$. Then the strict closure of $\, C$ in $M(A)=M$ is equal
to $C'\cap M\,$.
\endproclaim

{\bf Proof.} According to Lemma 1 (vi), we have to show that $C$ 
contains a two-sided approximate unit for $A\,$. Without loss of
generality we may assume that $A\ne \{0\}\,$, hence $C\ne\{0\}\,$.\par
Let $(e_\iota )_{\iota \in I}$ be a maximal family of mutually
orthogonal non-zero projections in $C\,$. Then
$$
  \bigvee_\iota e_\iota =1_M.
$$
Indeed, $e_o=1_M-\bigvee_\iota e_\iota $ belongs to $C'\cap M,$ so
Lemma 4 (iv) yields $e_o=\bigvee\{f\in C; f\le e_o\hbox{ projection}\}.$
Thus $e_o\ne 0$ would imply the existence of some projection $0\ne f\le e_o$
in $C,$ contradicting the maximality of $(e_\iota )_{\iota \in I}\, .$

Denoting by $Z$ the centre of $M,$ we call central partition of $1_M$ any set
of mutually orthogonal projections in $Z$ with least upper bound $1_M.$ 
The projections
\medskip
\centerline{$\displaystyle
\bigvee_{p\in\Cal P} \Big(\sum_{\iota \in I_p}e_\iota\Big) p\, ,\quad 
\Cal P \hbox{ central partition of } 1_M\, , \quad
I_p\subset I\hbox{ finite for any } p\in\Cal P$}
\smallskip

\noindent belong to $\, C'\cap M$ and are finite (see [Be],\S 15,
Proposition 8), hence they belong to $,\ C'\cap A=C.$ We show that
their family is an (increasing positive) approximate unit for
$A\, .$ For we have to prove that every finite projection $e$ in
$M$ has the property
\medskip

\noindent (P)\hskip1.885cm$\cases \,\hbox{for every } \ve >0
\hbox{ there are $\Cal P$ and $I_p\, , p\in\Cal P$, with}\cr
\,\Big\|\Big( 1_M-\bigvee\limits_{p\in\Cal P}\big(
\sum\limits_{\iota\in I_p}e_\iota\big) p\Big) e\Big\|\le\ve\, .
\endcases$
\smallskip

\noindent But standard arguments show that every finite projection
$e$ in $M$ is of the form
\medskip
\centerline{$\displaystyle e=\bigvee_{n\ge 1}(e_{n,1}+\dots
+e_{n,n}) p_n\, ,$}
\smallskip

\noindent where $p_n, n\ge 1$ are mutually orthogonal projections
in $Z$ and, for every $n\ge 1, e_{n,1},\dots, e_{n,n}$ are mutually
orthogonal abelian projections of central support $p_n$ (use
[Be], \S 18, Exercises 3, 4 and Proposition 1), so it is enough to
prove (P) for every abelian projection $e$ in $M.$ Moreover, since
every abelian projection is majorized by an abelian projection of
central  support $1_M,$ without loss of generality we can restrict
us to the case of an abelian projection $e$ of central support
$1_M\,$.

For every $x\in M$ there exists a unique $\Phi_e(x)\in Z$ such that
\medskip
\centerline{$exe=\Phi_e(x)e$}
\medskip

\noindent
(see [Be], \S 15, Proposition 6 and \S 5). Clearly, $\Phi_e: M\to Z$
is a conditional expectation and, according to [Kap 2], Lemma 7,
it is completely additive on the projection lattice of $M\,$.
Furthermore, $Z\ni z\mapsto ze\in Ze$ being $*$-isomorphism, we have
\medskip
\centerline{$\|xe\|^2=\|ex^*xe\|=\|\Phi_e(x^*x)e\|=
\|\Phi_e(x^*x)\|,\quad x\in M\, .$}
\medskip

Now, by the complete additivity of $\Phi_e,$
\medskip
\centerline{$\displaystyle
\sum_\iota  \Phi_e(e_\iota )=\Phi_e(1_M)=1_M\, .$}
\smallskip

\noindent
Thus, according to [Kap 2], Lemma 5, for every $e>0$ there
exist a central partition $\Cal P$ of $1_M$ and finite sets
$I_p\subset I, p\in \Cal P$ such that
\medskip
\centerline{$\displaystyle
\Big\|\Big( 1_M-\sum_{\iota \in I_p}\Phi_e(e_\iota )\Big) p\,\Big\|
\le\ve^2\hbox{ for all }p\in\Cal P\, .$}
\smallskip

\noindent
But then we have for every $p\in \Cal P$
$$\aligned
\Big\|\Big( 1_M-\sum_{\iota \in I_p}e_\iota\Big) pe\,\Big\|^2&=
\Big\|\Phi_e((1_M-\sum_{\iota \in I p}e_\iota\Big) p\,\Big\|=\\
&=\Big\|\Big( 1_M-\sum_{\iota \in I p}\Phi_e (e_\iota)\Big) p
\,\Big\|\le\ve^2,
\endaligned
$$
\noindent
so, taking into account [Kap 1], Lemma 2.5,
\medskip
\centerline{$\displaystyle
\Big\|\Big( 1_M-\bigvee_{p\in\Cal P}\big(\sum_{\iota \in I_p}
e_\iota\big) p\Big) e\,\Big\| = \sup_{p\in\Cal P}
\Big\|\Big( 1_M-\sum_{\iota \in I_p}e_\iota\Big) pe\,\Big\|
\le\ve\, .$}

\hfill$\square$
\bigskip

{\bf 2  Abelian Strict Closure in Continuous AW*-algebras }
\bigskip

For the treatment of the case of continuous $M$ we need several
lemmas on \AWs~, which could be of interest for themselves. First
we extend [Z], Lemma 2.2, concerning a Darboux property of normal
functionals on von Neumann algebras without minimal projections,
to the case of centre valued completely additive maps on the
projection lattice of a continuous \AW (similar results can be
found in [Ars-Z] and, for tracial maps, in [Kad 2], Prop. 3.13,
[Kaf], Prop. 27).

\medskip

\proclaim{Lemma 2 } Let $M$ be a continuous $AW^*$-algebra, $Z$ its
centre, $C$ a masa of $M\,$, and
\smallskip

\centerline{
$\Phi:\{e\in M; e\hbox{ projection }\}\to Z^+$}
\medskip

\noindent a completely additive map such that
\medskip

\centerline{
$\Phi (ep)=\Phi(e)p,\quad e\in M\hbox{ and } p\in
Z\hbox{ projections.}$}
\medskip

\noindent Then, for every projection $e\in C\,$,
\medskip

\centerline{
$\{z\in Z; 0\le z\le\Phi(e)\}=\{\Phi(f); e\ge f\in
C\hbox{ projection}\}\, .$}
\endproclaim
\medskip

{\bf Proof.} a) First we prove that for every projection
$0\ne g\in C$ there exists a projection $0\ne h\le g$ in $C$
such that
\medskip

\centerline{
$\Phi(h)\le\frac12\Phi(g)\, .$}
\medskip

\noindent The case $\Phi(g)=0$ being trivial, we can assume
without loss of generality that $\Phi(g)\ne 0.$\par

Let $(g_\iota )_\iota $ be a maximal family of mutually
orthogonal  projections in $Cg$ such that $\Phi(g_\iota )=0$
for every $\iota\, .$ Put $g_1=g-\bigvee_\iota g_\iota\in C\, .$
Then
\medskip

\centerline{
$\displaystyle \Phi(g_1)=\Phi(g)-\sum_\iota\Phi(g_\iota)=
\Phi(g)\ne 0\, ,$}
\smallskip

\noindent so $g_1\ne 0.$ By the maximality of $(g_\iota)_\iota,$
for no projection $0\ne g'\le g_1$ in $C$ can hold $\Phi(g')=0\,$.

Now there exists a projection $g_2\le g_1$ in $C$ such that
$g_2\notin Zg_1\,$. For let us assume the contrary, that is that
\medskip

\centerline{
$C=Zg_1+C(1_M-g_1)\, .$}

\noindent There exist projections $h_1,h_2\in M$ such that
$g_1=h_1+h_2$ and $h_1\sim h_2$ ([Be], \S 19, Th. 1) and then
\smallskip

\centerline{
$C\subset Zh_1+Zh_2+C(1_M-g_1)$}
\medskip

\noindent and the maximal abelianness of $C$ imply that
\medskip

\centerline{
$C=Zh_1+Zh_2+C(1_M-g_1)\, .$}

\noindent Thus
\smallskip

\centerline{
$h_1,h_2\in C g_1=Zg_1\, .$}
\medskip

\noindent But, denoting by $z(g_1)$ the central support of $g_1,$
\medskip

\centerline{
$Z z(g_1)\ni z\mapsto zg_1\in Zg_1\, .$}
\medskip

\noindent is a $*$-isomorphisms and it follows that $h_1$ and
$h_2$ have orthogonal central supports, in contradiction to
$h_1\sim h_2\ne 0\,$.

We claim that $\Phi (g_2)\Phi (g_1-g_2)\ne 0.$ Indeed, otherwise
it would exist a projection $p\in Z$ such that 
\medskip

\centerline{
$\Phi (g_2)=\Phi (g_2)p\hbox{ and } \Phi (g_1-g_2)p=0$}
\medskip

\noindent and it would follow successively
$$\aligned
\Phi (g_2(1_M-p))&=0\hbox{ and } \Phi ((g_1-g_2)p)=0\, ,\\
(g_2(1_M-p))&=0\hbox{ and } (g_1-g_2)p=0\, ,\\
g_2=g_2p&=g_1p\in Zg_1\, .
\endaligned
$$

Let $q\in Z$ denote the support projection of
$(\Phi (g_1)-2\Phi (g_2))_+.$ Then
$$\aligned
\Phi (g_1q)-2\Phi (g_2q)\,&=(\Phi (g_1)-2\Phi (g_2))_+\ge 0,\\
\Phi (g_2q)\le\frac12\Phi (g_1q)&\le\frac12\Phi (g_1)\le
\frac12\Phi (g)\, .
\endaligned
$$
\noindent
Similarly,
$$
\Phi ((g_1-g_2)(1_M-q))\le\frac12\Phi (g)\, .
$$
But we can not have simultaneously
$$
\Phi(g_2q)=0\hbox{ and } \Phi((g_1-g_2)(1_M-q))=0\, ,
$$
\noindent
because this would imply
$$
\Phi(g_2)\Phi(g_1-g_2)=
\Phi(g_2q)\Phi(g_1-g_2)+\Phi(g_2)\Phi((1_M-q)(g_1-g_2))=0\, .
$$
Therefore, putting $h=g_2q$ if $\Phi(g_2q)\ne 0$ and
$h=(g_1-g_2)(1_M-q)$ otherwise, $h$ is a non-zero projection in
$C\,$, majorized by $g\,$, such that $\Phi(h)\le\frac12\Phi(g)
\, .$
\smallskip

b) Now let $e\in C$ be a projection and let $x\in Z\, ,\, 0\le
z\le\Phi(e)$ be arbitrary. Choose a maximal family $(f_\iota)_\iota$
of mutually orthogonal projections in $Ce$ satisfying
\medskip

\centerline{
$\displaystyle \sum_\iota\Phi(f_\iota)\le z\, .$}
\smallskip

\noindent Then the projection $f=\bigvee_\iota f_\iota\le e$
belongs to $C$ and
\medskip

\centerline{
$\displaystyle \Phi(f)=\sum_\iota\Phi(f_\iota)\le z\, .$}
\smallskip

\noindent We claim that actually $\Phi(f)=z\, .$

For let us assume the contrary. Then there exist a projection
$0\ne p\in Z$ and $\ve>0$ such that
\medskip

\centerline{
$(z-\Phi(f))p\ge\ve p\, .$}
\medskip

\noindent The projection $g=(e-f)p\in C$ is not zero, because
otherwise it would follow
\medskip

\centerline{
$0=(\Phi(e)-\Phi(f))p\ge (z-\Phi(f))p\ge\ve p\, ,$}
\medskip

\noindent contradicting $p\ne 0,\ve >0\, .$ Choosing an
integer $n\ge 1$ with $2^{-n}\|\Phi(e-f)\|\le \ve\, ,\, n$-fold
application of a) yields the existence of a projection
$0\ne h\le g$ in $C$ such that
\smallskip

\centerline{
$\Phi(h)\le 2^{-n}\Phi((e-f)p)\le\ve p\, .$}
\medskip

\noindent Since $0\ne h\in Ce$ is orthogonal to every
$f_\iota$ and
\medskip

\noindent\hskip2.54cm$\displaystyle \Phi(h)+\sum_\iota
\Phi (f_\iota)=\Phi (h)+\Phi (f)\le\ve p+\Phi(f)\le z\, ,$
\smallskip

\noindent the maximality of $(f_\iota)_\iota$ is contradicted.

\hfill$\square$
\bigskip

It is well known that if the projection family $(e_\iota)_\iota$
in a finite \AW $M$ is upward directed and, for some projection
$f\in M\, ,\, e_\iota\prec f$ for all $\iota\, ,$ then
$\bigvee_\iota e_\iota\prec f$ (see [Be], \S 33, Exercise 1).
The above statement actually holds in any \AW $M$ under the only
assumption of the finiteness of $f$ (see Appendix, Cor. 1).
Here we give a proof for this, assuming additionally that the
projections $e_\iota$ are the finite partial sums of a family
of mutually orthogonal projections in $M\, :$
\medskip

\proclaim{Lemma 3} Let $M$ be an $AW^*$-algebra, $f\in M$ a
finite projection, and $(e_\iota)_{\iota\in I}$ a family of
mutually orthogonal projections in $M$ such that
\medskip

\centerline{
$\displaystyle \sum_{\iota\in F}e_\iota\prec
f\hbox{ for every finite } F\subset I\, .$}
\smallskip

\noindent Then
\smallskip

\centerline{
$\displaystyle \bigvee_{\iota\in I} e_\iota\prec f\, .$}
\endproclaim

{\bf Proof} According to the theory of Murray-von Neumann
equivalence for projections in $AW^*$-algebras, we can assume
without loss of generality that either $fMf$ is of type $I_n$
for some natural number $n\ge 1\, ,$ or that it is continuous
(see [Be], \S 15, Th.1, \S 18, Th. 2, \S 6, Cor. 2 of Prop. 4).\par

Let us first assume that $fMf$ is of type $I_n.$ By the Zorn Lemma
there exists a maximal set $\Cal P$ of mutually orthogonal
central projections in $M$ such that 
\medskip

\centerline{
$\card \{\iota\in I; pe_\iota\ne 0\}\le n\hbox{ for every } p
\in \P\, .$}
\medskip

We claim that $\bigvee\P=1_M\, .$ For let us assume that
$p_o=\bigvee\Cal P\ne 1_M\, .$ Then we can find recursively $n+1$
indices $\iota_1,\dots,\iota_{n+1}\in I$ such that
\medskip

\centerline{
$p_1=(1_M-p_o)z(e_{\iota_1})\dots z(e_{\iota_{n+1}})\ne 0\, ,$}
\medskip

\noindent where $z(e_\iota)$ denotes the central support of
$e_\iota\, .$ By the assumption of the lemma there exist mutually
orthogonal projections $f_{\iota_1},\dots ,f_{\iota_{n+1}}\le f$
in $M$ such that $e_{\iota_j}\sim f_{\iota_j}$ for every
$1\le j\le n+1\, .$ For every $1\le j\le n+1\,$, the central
support of $p_1f_{\iota_j}$ is $p_1\,$, so there exists an
abelian projection $g_j\le p_1 f_{\iota_j}$ of central support
$p_1$ (see [Be], \S 18, exercise 4). But then $g_1, \dots, g_{n+1}$
are mutually orthogonal, equivalent, non-zero projections in $fMf$
(see [Be], \S 18, Prop.1), which contradicts [Be], \S 18, Prop. 4.

By the very orthogonal additivity of equivalence in \AWs (see [Be], 
\S 11, Prop. 2) we conclude that
\medskip

\centerline{
$\displaystyle \bigvee_{\iota\in I}e_\iota=
\bigvee\bigg\{\sum_{pe_\iota\ne 0} pe_{\iota_j}\, ;\, p\in
\Cal P\bigg\}\prec\bigvee\{pf\, ;\, p\in \Cal P\}=f\, .$}
\medskip

Let us next assume that $fMf$ is continuous and let
$x\mapsto x^\natural$ denote the centre valued dimension
function of the finite \AW $fMf$ (see [Be], Ch.6).

For every $\iota\in I$ there exists a projection $e'_\iota\le f$
in $M$ such that $e_\iota\sim e'_\iota\, .$ Since
$(e'_\iota)^\natural$ does not depend on the  choice of
$e'_\iota\,$, we can put
\medskip

\centerline{
$e_\iota^\natural=(e'_\iota)^\natural\, .$}
\medskip

By the assumption of the lemma, for every finite $F\subset I$
we can choose the projections $e'_\iota\, ,\,\iota\in F\,$,
mutually orthogonal and then
\medskip

\centerline{
$\displaystyle \sum_{\iota\in F}e_\iota^\natural=
\sum_{\iota\in F} (e'_\iota)^\natural=\bigg(\sum_{\iota\in F}
e'_\iota\bigg)^\natural\le f\, .$}
\medskip

\noindent It follows that all sums
\smallskip

\centerline{
$\displaystyle \sum_{\iota\in J}e_\iota^\natural\le f,\quad
J\subset I$}
\medskip

\noindent exist in the monotone complete centre of $fMf\, .$

Now let us consider the set of all families of mutually
orthogonal projections in $fMf$
$$
(f_\iota)_{\iota\in J} \hbox{ with } J\subset I\, ,
$$
for which $f_\iota\sim e_\iota$ for every $\iota\in J\,$.
We can endowe this set with the partial order
\medskip

\centerline{
$(f_\iota)_{\iota\in J}\le(f'_\iota)_{\iota\in J'}\,
\Longleftrightarrow\, J\subset J'\hbox{ and }f_\iota=f'_\iota
\hbox{ for all }\iota\in J\, .$}
\medskip

\noindent By the Zorn lemma there exists a maximal element
$(f_\iota)_{\iota\in J}$ of the above partially ordered set.
We claim that then $J=I\, .$ For let us assume the existence
of some $\iota_o\in I\backslash J\, .$ Since
\smallskip

\centerline{
$\displaystyle e_{\iota_o}^\natural +\bigg(\bigvee_{\iota\in J}
f_\iota\bigg)^\natural =e_{\iota_o}^\natural +\sum_{\iota\in J}
f_\iota^\natural\le\sum_{\iota\in I}e_\iota^\natural\le f\, ,$}
\medskip

\noindent that is
\medskip

\centerline{
$\displaystyle e_{\iota_o}^\natural\le\bigg( f-\bigvee_{\iota\in J}
f_\iota\bigg)^\natural\, ,$}
\medskip

\noindent by [Be], \S 33, Th.3 (particular case of the above
Lemma 5) there exists a projection $f_{\iota_o}\le
f-\bigvee_{\iota\in J}f_\iota$ in $M$ such that
$f_{\iota_o}^\natural =e_{\iota_o}^\natural =\left(e'_{\iota_o}
\right)^\natural\,$, hence $f_{\iota_o}^{}\sim e'_{\iota_o}
\sim e_{\iota_o}^{}\, .$ But this contradicts the maximality of
$(f_\iota)_{\iota\in J}.$

By the general additivity of equivalence in \AWs
(see [Be], \S 20, Th. 1) we can conclude also in this case that
\medskip

\centerline{
$\displaystyle \bigvee_{\iota\in I} e_\iota\sim\bigvee_{\iota\in I}
f_\iota\le f\, .$}

\hfill$\square$
\bigskip

Let $M$ be a semi-finite $AW^*$-algebra, and $A$ the norm-closed
linear span of all finite projections of $M\, .$ We recall that
then $M=M(A)\, .$\par

Let us call a masa $\tilde C$ of $M \ M$-semi-finite if
$\tilde C\cap A$ is an essential ideal of $\tilde C$ or, equivalently,
if every non-zero projections  in $\tilde C$ majorizes a non-zero
projection in $\tilde C\cap A$ (cf. with [Kaf]. Def. 1). For
$\tilde C\subset M$ are equivalent:\par 
\quad 1) $\tilde C$ is an $M$-semi-finite masa of $M\, ;$\par
\quad 2) $\tilde C=C'\cap M$ for some masa $C$ of $A\, .$\par 
\noindent 
Indeed, 2) implies 1) by Lemma 1 (iv), while 1) $\Rightarrow$ 2)
follows by noticing that, according to the $M$-semifiniteness of
$\tilde C\,$, every projection in $\tilde C$ is the least upper
bound of a family of mutually orthogonal projections from
$C=\tilde C\cap A\,$, and so $C'\cap M=\tilde C'\cap M=\tilde C
\, ,\, C'\cap A=(C'\cap M)\cap A=\tilde C\cap A=C\, .$\par

The following result extends [Kad 2], Th. 3.18 and [Kaf], Cor. 31
in the case of an $M$-semifinite masa :
\medskip

\proclaim{Theorem 2 {\rm $($on labeling Murray-von Neumann
equivalence classes$)$}} Let $M$ be a semi-finite $AW^*$-algebra,
$A$ the norm closed linear span of all finite projections of $M\,$,
and $C$ a masa of $A\,$. Then

\item{{\rm (i)}} for any projections $M\ni f\le e\in C'\cap M$
there exists a projection $f\sim g\le e$ in $C'\cap M\, ;$\par

\item{{\rm (ii)}} for any projections $M\ni f\le e\in C'\cap M$
of equal central supports, $f$ finite and $e$ properly infinite,
there is a set $\Cal P$ of mutually orthogonal central projections
in $M$ with $\bigvee\Cal P=1_M$ such that, for every $p\in \Cal P
\, ,\, ep$ is the least upper bound in the projection lattice of
$M$ of some family of mutually orthogonal projections from $C\,$,
each one of which is equivalent in $M$ to $fp\,$.\par
\endproclaim

{\bf Proof.} (a) First we prove (i) in the case $e\in C\,$.
Similarly as in the proof of Lemma 3, we can assume without
loss of generality that either $eMe=eAe$ is of type $I_n$ for
some natural number $n\ge 1\,$, or it is continuous.

If $eMe$ is of type $I_n\,$, by [Kad 2], Lemma 3.7 there exist
mutually orthogonal projections $e_1,\dots, e_n\in C$ with
$\sum\limits^n_{j=1}e_j=e\,$, such that each $e_j$ is abelian
in $M$ and has the same central support in $M$ as $e$ (actually
[Kad 2], Lemma 3.7 is proved only for von Neumann algebras, but
an inspection of the proof shows that it works without any change
also in the realm of the $AW^*$-algebras). On the other hand,
using [Be], \S 18, Exercise 4 and Prop. 4, it is easy to see that
there exist mutually orthogonal abelian projections
$f_1,\dots,f_n\in M$ with $\sum\limits^n_{j=1}f_j=f$ and central
supports $z(f)=z(f_1)\ge\dots\ge z(f_n)\,$.
By [Be], \S 18, Prop. 1 it follows that $f_j\sim e_j z(f_j)$
for all $1\le j\le n\,$, so $f$ is equivalent to
$C\ni \sum\limits^n_{j=1} e_j z(f_j)\le e\,$.\par

Now let us assume that $eMe$ is continuous and let
$x\mapsto x^\natural$ denote the centre valued dimension function
of finite \AW $eMe\,$. Then Lemma 2 yields the existence of a
projection $C\ni g\le e$ such that $g^\natural=f^\natural\,$,
hence $g\sim f\,$.
\smallskip

(b) Next we prove (i) in the case $f\in A\,$.

By Lemma 1 (iv) there exists a family $(e_\iota)_{\iota\in I}$
of mutually orthogonal projections in $C$ such that
\smallskip

\centerline{
$\displaystyle e= \bigvee_{\iota\in I}e_\iota\, .$}
\medskip

\noindent Let $\Cal P$ be a maximal set of mutually orthogonal
central projections in $M$ such that, for every $p\in \Cal P$,
there is a finite set $F_p\subset I$ with
\medskip

\centerline{
$\displaystyle fp\prec p\sum_{\iota\in F_p}e_\iota\in C\, .$}
\smallskip

\noindent By the above part (a) of the proof, for every $p\in\P$
there exists a projection
\smallskip

\centerline{$g(p)\in C\hbox{ with }fp\sim g(p)\le p
\sum\limits_{\iota\in F_p}e_\iota\, .$}
\noindent If $\,\bigvee\Cal P=1_M$ then $f=\bigvee\{fp\,;\, p
\in \Cal P\}$ is equivalent to $C'\cap M\ni\bigvee\{g(p)\, ;\,
p\in\Cal P\}\le e\,$, so let us assume in the sequel that
$p_o=1_M-\bigvee\Cal P\ne 0\,$.

By the maximality of $\Cal P$ and by the comparison theorem
(see [Be], \S 14, Cor. 1 of Prop. 7) we have
\smallskip

\centerline{
$\displaystyle p_o\sum_{\iota\in F}e_\iota \prec
f \hbox{ for every finite } F\subset I\, .$}
\medskip

\noindent According to Lemma 3 it follows that
\medskip

\centerline{
$\displaystyle p_o e=\bigvee_{\iota\in I}p_oe_\iota\prec f\, ,$}
\smallskip

\noindent so by the Schr\"oder-Bernstein theorem (see [Be],\S 12)
we have
\medskip

\centerline{
$fp_o\sim ep_o\, .$}
\medskip

\noindent Consequently $f=fp_o +\bigvee\{fp\, ;\, p\in\Cal P\}$
is equivalent to
\medskip

\centerline{
$C'\cap M\ni ep_o +\bigvee\{g(p)\, ;\, p\in \Cal P\}\le e\,$.}
\medskip

(c) Now we prove (ii).

Let $\P$ be a maximal set of mutually orthogonal central
projections in $M$ such that, for every $p\in\P\, ,\, ep$ is
the least upper bound in the projection lattice of $M$ of some
family of mutually orthogonal projections from $C\,$, each one
of which is equivalent in $M$ to $fp\, .$ We claim that then
$\bigvee\P=1_M\, .$

For let us assume that $p_o=1_M-\bigvee\P\ne 0\, .$ We notice
that $fp\ne 0$ for any central projection $0\ne p\le p_o$ in
$M\, :$ indeed, otherwise $p$ would be orthogonal to the common
central support of $f$ and $e\,$, so $ep=0$ would be equal to
$fp=0\in C\,$, in contradiction with the maximality  of $\P\,$.

Let $(e_\iota)_{\iota\in I}$ be a maximal family of mutually
orthogonal projections in $C$ such that $fp_o\sim e_\iota\le ep_o$
for all $\iota\in I.$
By the comparison theorem there exists a central projection
$p_1\le p_o$ in $M$ such that
\medskip
\centerline{$\displaystyle \Big( ep_o-\bigvee_{\iota\in I}
e_\iota\Big) p_1\prec fp_1\, ,$}
\smallskip

\centerline{$\displaystyle \Big( ep_o-\bigvee_{\iota\in I}
e_\iota\Big) (p_o-p_1)\succ f(p_o-p_1)\, .$}
\medskip

\noindent Then $p_1\ne 0\, :$ indeed, $p_1=0$ would imply
$$
A\ni fp_o\prec ep_o-\bigvee_{\iota\in I}e_\iota\in C'\cap M
$$
and, by the above proved (b), it would exist a projection
$fp_o\sim e'\le ep_o -\bigvee_{\iota\in I}e_\iota$ in
$(C'\cap M)\cap A=C\,$, contradicting the maximality of
$(e_\iota)_{\iota\in I}\,$. Put
$$
e_o=ep_1-\bigvee_{\iota\in I}e_\iota p_1\prec fp_1\, .
$$
Then $e_o$ is finite and belongs to $C'\cap M\,$, so it
belongs to $C'\cap A=C\,$. On the other hand, the proper
infiniteness of $e$ and $ep_1\ne 0$ imply that
$ep_1=e_o+\bigvee_{\iota\in I}e_\iota p_1$ is properly infinite.
It follows that the set $I$ is necessarily infinite, hence
containing an infinite  sequence $\iota_1,\iota_2,\,\dots\;$.

For every $j\ge 1\, ,\, e_o\prec fp_1\sim e_{\iota_j}p_1\in C$
and the above proved a) yield the existence of some projection
$e_o\sim e^{(1)}_{\iota_j}\le e_{\iota_j}p_1$ in $C\,$.
In particular, all projections $e^{(1)}_{\iota_j}$ are
equivalent, hence, the projections $e_\iota p_1$ being finite,
the projections $e^{(2)}_{\iota_j}=e_{\iota_j}p_1-e^{(1)}_{\iota_j}$
are also all equivalent (see [Be], \S 17, Exercise 3). Consequently,
the projections from $C$
$$
e'_{\iota_1}=e_o+e^{(2)}_{\iota_1}\hbox{ and } e'_{\iota_j}=
e^{(1)}_{\iota_{j-1}}+e^{(2)}_{\iota_j}\, ,\quad j\ge 2
$$
are all equivalent in $M$ to $e^{(1)}_{\iota_1}+e^{(2)}_{\iota_1}=e_{\iota_1}
p_1\sim fp_1\,$. Clearly, they are mutually orthogonal and
$$
\bigvee_{j\ge 1}e'_{\iota_j}=e_o\lor\bigvee_{j\ge 1}
e^{(1)}_{\iota_j}\lor\bigvee_{j\ge 1}e^{(2)}_{\iota_j}=
e_o\lor\bigvee_{j\ge 1}e_{\iota_j}p_1\, .
$$
Letting 
$$
e'_\iota =e_\iota p_1\hbox{ for }\iota\in I\,\backslash\,
\{\iota_1,\iota_2,\,\dots\},
$$
we conclude that all projections $e_\iota'\, ,\,\iota\in I\,$,
belong to $C$ and are equivalent in $M$ to $fp_1\,$. Moreover,
they are mutually orthogonal and
$$
\bigvee_{\iota\in I}e'_\iota=\bigvee_{j\ge 1}e'_{\iota_j}
\lor\bigvee_{\iota\ne\iota_j}e'_\iota=e_o\lor\bigvee_{j\ge 1}
e_{\iota_j} p_1\lor\bigvee_{\iota\ne\iota_j}e_\iota p_1=
e_o\lor\bigvee_{\iota\in I}e_\iota p_1=ep_1\, .
$$
But this contradicts the maximality of $\P\,$.
\smallskip

(d) Finally we prove (i) in full generality.

We can assume without loss of generality that either $f$ is finite,
or it is properly infinite. The case of finite $f$ was already
settled in (b), so it remains to consider only the case of properly
infinite $f\,$.

Choose some finite projection $M\ni f_o\le f$ of the same
central support as $f$ (see [Be], \S 17, Exercise 19 iii)).
According to the above proved (c), we  can assume without loss of
generality that there are families  $(e_\iota)_{\iota\in I}$ and
$(f_\kappa )_{\kappa\in K}$ of mutually orthogonal projections in
$M$ such that
$$
e_\iota\sim f_o\sim f_\kappa\,\hbox{ for all }\,\iota\in I\,
\hbox{ and }\,\kappa\in K
$$
$$
\bigvee_{\iota\in I} e_\iota =e,\quad \bigvee_{\kappa\in K}f_\kappa
=f\, .
$$
If card $K\le$ card $I\,$, that is if there exists an injective
map $K\ni\kappa\mapsto\iota(\kappa )\in I\,$, then the projection
$g=\bigvee_{\kappa\in K}e_{\iota(\kappa )}\le e$ belongs to
$C'\cap M$ and is equivalent to $\bigvee_{\kappa\in K} f_\kappa
=f\,$. On the other hand, if card $I\le$ card $K\,$, then
$e=\bigvee_{\iota\in I} e_\iota\prec\bigvee_{\kappa\in K}f_\kappa
=f\le e$ and the Schr\"oder-Bernstein theorem imply that
$e\sim f\,$.

\hfill$\square$
\bigskip

Let us now prove the statement of [Kad 2], Th. 3.18 and [Kaf],
Cor. 31 in the case of an $M$-semifinite masa of an arbitrary
semifinite \AW $M\,$:
\medskip

\proclaim{Corollary} Let $M$ be a semifinite $AW^*$-algebra, $A$
the norm-closed linear span of all finite projections of $M\,$,
and $C$ a masa of $A\,$. If $e\in C'\cap M$ is a projection and
$1\le n\le \aleph_o$ is a cardinal number such that $e$ is the
least upper bound of $n$ mutually orthogonal, equivalent
projections from $M\,$, then there exist $n$ mutually orthogonal
projections in $C'\cap M\,$, all equivalent in $M\,$, whose least
upper bound is $e\,$.
\endproclaim

{\bf Proof.} It is enough to treat separately the case of finite
respectively properly infinite $e\,$. If $e$ is finite, $n$ can
be only a natural number. Let $f_1,\dots f_n$ be mutually orthogonal,
equivalent projections in $M$ with $\sum\limits_{j=1}^{n}f_j=e\,$.
By (i) in the above theorem there exists a projection
$f_1\sim e_1\le e$ in $C\,$. Since $e$ is finite, it follows that
$\sum\limits^n_{j=2}f_j\sim e-e_1\,$, so we can apply again (i)
in the above theorem to get a projection $f_2\sim e_2\le e-e_1$
in $C\,$. By induction we obtain $n$ mutually orthogonal
projections $e_1,\,\dots\, ,e_n\in C$ such that $f_j\sim e_j$
for all $j$ and $\sum\limits^n_{j=1}e_j=e\,$.

Now let us assume that $e$ is properly infinite and consider
a set $I$ of cardinality $n\,$. Choosing a finite projection
$M\ni f\le e$ of the same central support as $e$ (see [Be], \S 17,
Exercise 19 iii)), (ii) in the above theorem entails the existence
of a set $\P$ of mutually orthogonal central projections in $M$
with $\bigvee \P=1_M$ such that, for every $p\in \P\, ,\, ep$
is the least upper bound of some set $\Cal E_p$ of mutually
orthogonal projections from $C\,$, each one of which is
equivalent in $M$ to $fp\,$. If $ep\ne 0$ then $\Cal E_p$ must
be infinite, so there exists a partition
$(\Cal E_{p,\iota})_{\iota\in I}$ of $\Cal E_p$ in $n$ sets
of equal cardinality. Then the projections $e_\iota =
\bigvee_{ep\ne 0}\bigvee\Cal E_{p,\iota}\, ,\,\iota\in I\,$,
belong to $C'\cap M\,$, are mutually orthogonal and equivalent
in $M\,$, and $\bigvee_{\iota\in I}e_{\iota}=e\,$.

\hfill$\square$
\bigskip 

Let $M$ be a finite \AW with centre $Z$ and let $x\mapsto x^\natural$
denote its centre valued dimension function (see [Be], Ch. 6).
It is known (see [Bl-Ha], II, 1) that $\natural$ can be uniquely
extended to a centre valued quasitrace on $M,$ that is to a map
$\Phi:M\to Z$ such that

\item{-} $\Phi$ is linear on commutative $*$-subalgebras of $M\,$,

\item{-} $\Phi(a+ib)=\Phi(a)+i\Phi(b)$ for all selfadjoint $a, b\in 
M\,$,

\item{-} $\Phi$ acts identically on $Z\,$,

\item{-} $0\leqslant\Phi(x^*x)=\Phi(xx^*)$ for all $x\in M\,$,

\noindent and then

\item{-} $\Phi(a)\leqslant \Phi(b)$ whenever $a\leqslant b$ are
selfadjoint elements of $M\,$,

\item{-} $\Phi$  is norm continuous, more precisely, $\|\Phi(a)-\Phi(b)\|
\leqslant\|a-b\|$ for all selfadjoint $a, b\in M\,$.

\noindent We shall use the symbol $\natural$ to denote also
the above $\Phi\,$.

According to classical results of F.J. Murray and J. von Neumann,
the centre valued quasitrace of every finite \Wag is additive,
hence linear.

It is an open question, raised by I. Kaplansky, whether the
centre valued  quasitrace of every finite \AW is additive.
Recently U. Haagerup has proven that the answer to Kaplansky's
question is positive for any finite \AW which is generated
(as an $AW^*$-algebra) by an exact \Csalg (see [Haa], Th. 5.11, 
Prop. 3.12 and Lemma 3.7 (4)).

We notice that if $M$ is a finite \AW and $n\geq 1$ is an integer,
then the $*$-algebra Mat$_n(M)$ of all $n\times n$ matrices over
$M$ is again a finite \AW (see [Be], \S 62). Denoting by
$\natural$ and $\natural_n$ the respective centre valued
quasitraces, it is easily seen that
\smallskip

\centerline{$\displaystyle n\cdot
\pmatrix x&0&&0\cr 0&0&&\cr &&\ddots&\cr 0&&&0\endpmatrix^{\natural_n}
=\pmatrix x^\natural&0&&0\cr 0&x^\natural&&\cr &&\ddots&\cr 0&&&
x^\natural\endpmatrix\, , \quad x\in M\, .$}
\medskip

\noindent Moreover the additivity of $\natural$ is equivalent
with the validity of
\medskip

\centerline{$\displaystyle 2\cdot\pmatrix x_{11}& x_{12}\cr
x_{21}& x_{22}\endpmatrix^{\natural_2}=
\pmatrix x_{11}^\natural +x_{22}^\natural&0\cr
0&x_{11}^\natural+x_{22}^\natural\endpmatrix\, .$}
\medskip

\noindent Indeed, using the above equality, we get for all
$0\le a, b\in M$
$$
\aligned
\pmatrix(a+b)^\natural&0\cr 0&(a+b)^\natural\endpmatrix
&=2\cdot\pmatrix  a+b&0\cr0&0\endpmatrix^{\natural_2}=\\
&=2\cdot\left[\pmatrix a^{1/2}&b^{1/2}\cr0&0\endpmatrix
\pmatrix a^{1/2}&0\cr b^{1/2}&0\endpmatrix\right]^{\natural_2}=\\
&=2\cdot\left[\pmatrix a^{1/2}&0\cr b^{1/2}&0\endpmatrix 
\pmatrix a^{1/2}&b^{1/2}\cr0&0\endpmatrix\right]^{\natural_2}=\\
&=2\cdot\pmatrix a&a^{1/2}b^{1/2}\cr b^{1/2}a^{1/2}&b
\endpmatrix^{\natural_2}=\\
&=\pmatrix a^\natural+b^\natural&0\cr 0&a^\natural+b^\natural
\endpmatrix.
\endaligned
$$
Conversely, assuming that $\natural$ is additive, it is
easy to verify that
\medskip

\centerline{$\displaystyle \pmatrix x_{11}&x_{12}\cr x_{21}&x_{22}
\endpmatrix\mapsto\frac12 \pmatrix x_{11}^\natural+x^\natural_{22}
&0\cr 0&x^\natural_{11}+x^\natural_{22}\endpmatrix$}
\smallskip

\noindent is a centre  valued quasitrace on Mat$_2(M)\,$.

For a given $\delta>0$, we say that the centre valued
quasitrace $\natural$ of a finite \AW $M$ is $\delta$-subadditive
(resp. $\delta$-superadditive) if the map $M_+\ni a\mapsto
(a^\natural)^\delta$ is subadditive (resp. superadditive).
Clearly, $\delta$-subadditivity ($\delta$-superadditivity) of
$\natural$ implies its $\delta'$-subadditivity
($\delta'$-superaddivity) whenever $\delta'<\delta (\delta'>
\delta).$ It was proven by U. Haagerup that $\natural$ is always
$\frac12$-subadditive (see [Haa], Lemma 3.5 (1)) and it seems
reasonable to conjecture that it is also always 2-superadditive
(or, at least, $k$-superadditive for some $k\geqslant 1$).

We notice as a curiosity that, for any two projections $p\, ,\, q$
in a finite \AW $M$ with centre valued quasitrace $\natural\,$,
\smallskip

\centerline{$(p+q)^\natural=p^\natural +q^\natural\, .$}
\smallskip

\noindent Indeed, since
\smallskip

\centerline{$\pmatrix p+q&0\cr 0&0\endpmatrix =\pmatrix
p&\pm q\cr 0&0\endpmatrix\pmatrix p&0\cr\pm q&0\endpmatrix\, ,$}
\smallskip

\centerline{$\pmatrix p&\pm pq\cr \pm qp& q\endpmatrix=
\pmatrix p&0\cr \pm q&0\endpmatrix \pmatrix p&\pm q\cr
0&0\endpmatrix\, ,$}
\smallskip

\noindent and $\pmatrix p&pq\cr qp&q\endpmatrix,\pmatrix p
&-pq\cr -qp&q\endpmatrix$ commute, we have
\medskip

\noindent\hskip1.828cm
$\pmatrix (p+q)^\natural&0\cr 0&(p+q)^\natural\endpmatrix =
2\pmatrix p+q&0\cr 0&0\endpmatrix^{\natural_2}=$
\smallskip

\noindent\hskip5.275cm
$=\pmatrix p&pq\cr qp&q\endpmatrix^{\natural_2}+
\pmatrix p&-pq\cr -qp&q\endpmatrix^{\natural_2}=$
\smallskip

\noindent\hskip5.275cm
$=2\pmatrix p&0\cr 0&q\endpmatrix^{\natural_2}=$
\smallskip

\noindent\hskip5.275cm
$=2\pmatrix p&0\cr 0&0\endpmatrix^{\natural_2}+2
\pmatrix 0&0\cr 0&q\endpmatrix^{\natural_2}=$
\smallskip

\noindent\hskip5.275cm
$=\pmatrix p^\natural&0\cr 0&p^\natural\endpmatrix+
\pmatrix q^\natural&0\cr 0& q^\natural\endpmatrix=$
\smallskip

\noindent\hskip5.275cm
$=\pmatrix p^\natural+q^\natural&0\cr 0&p^\natural+
q^\natural\endpmatrix\, .$
\medskip

\noindent This can be deduced also from Haagerup's result,
taking to account that the \Cag generated by two projections
is of type $I\,$, hence nuclear, hence exact.
\medskip

\proclaim{Lemma 4} Let $M$ be a finite $AW^*$-algebra, whose
centre valued quasitrace $\natural$ is $k$-superadditive for
some $k\ge 1\, .$ Let further $e_1\, ,\,\dots\, ,\, e_n\in M$
be mutually equivalent projections with $\sum\limits^n_{j=1}
e_j=1_M\,$. Then there exists a projection $e_1\sim p\in M$
such that, for every projection $f\in
\{e_1\, ,\,\dots\, ,\, e_n\}'\cap M\,$,
$$
f^\natural \ge (1-\|(1_M-f)p\|^2)n^{\frac1k-1}1_M\, .
$$
\endproclaim
\medskip

{\bf Proof.} Let $v_1,\dots v_n\in M$ be partial isometries
such that
$$
v^*_j v_j=e_1,\quad v_jv^*_j=e_j,\quad 1\le j\le n\, .
$$
Since
$$
\bigg(\frac1{\sqrt n}\sum^n_{j=1}v_j\bigg)^*
\frac1{\sqrt n}\sum^n_{j=1}v_j=\frac1n\sum^n_{j_1,j_2=1}
v^*_{j_1}v_{j_2}=\frac1n \sum^n_{j=1}v^*_jv_j=e_1\, ,
$$
$\displaystyle p=\frac1n\!\sum^n_{j_1,j_2=1}\!\!v_{j_1}v^*_{j_2}=
\frac1{\sqrt n}\sum^n_{j=1}v_j\bigg(\frac1{\sqrt n}
\sum^n_{j=1}v_j\bigg)^*\!$ is a projection in $M$ equivalent
to $e_1\,$.
\smallskip

Now let the projection
\smallskip

\centerline{$f\in\{e_1\, ,\,\dots\, ,\, e_n\}'\cap M$}
\medskip

\noindent be arbitrary and set $\delta=\|(1_M-f)p\|\,$. Since
the case $\delta=1$ is trivial, we can assume without loss of
generality that $\delta <1\,$. Then 
\medskip

\centerline{$\|p-pfp\|=\|(1_M-f)p\|^2=\delta^2<1\, ,$}
\medskip

\noindent so $pfp\geq (1-\delta^2)p$ is invertible in $pMp\,$.
Thus the polar  decomposition $fp=w\cdot|fp|$ exists in the
\Cag generated by $p$ and $f$ and we have
\medskip

\centerline{$w^*w=p, fpf=w(pfp)w^*\geq (1-\delta^2)ww^*\, .$}
\medskip

\noindent Let us denote $\zeta=e^{i\frac n\pi}\, .$ Then
\medskip

\centerline{$u=\sum^n_{j=1}\zeta^j e_j\in
\{e_1,\dots,e_n, f\}'\cap M$}
\medskip

\noindent is unitary. Since

$$\aligned
u^mpu^{-m}&=\frac1n\sum^n_{j, j_1,j_2,j'=1}\zeta^{mj}e_jv_{j_1}
v^*_{j_2}\zeta^{-mj'}e_{j'}=\\
&=\frac1n\sum^n_{j_1, j_2=1}\zeta^{m(j_1-j_2)}v_{j_1}v^*_{j_2}
\endaligned
$$
and
\medskip

\centerline{$\sum^n_{m=1}\zeta^{mj}=0\hbox{ for every }
1\le j\le n-1\, ,$}
\medskip

\noindent we have
$$
\aligned
\sum^n_{m=1}u^mpu^{-m}&=\frac1n\sum^n_{j_1,j_2=1}\bigg(\sum^n_{m=1}
\zeta^{m(j_1-j_2)}\bigg)v_{j_1}v^*_{j_2}=\\
&=\frac1n\sum^n_{j_1=1}nv_{j_1}v^*_{j_1}=1_M
\endaligned
$$
Therefore
\medskip

\centerline{$\displaystyle f=f\sum^n_{m=1}u^mpu^{-m}f=\sum^n_{m=1}
u^m(fpf)u^{-m}\ge(1-\delta^2)\sum^n_{m=1}u^mww^*u^{-m}$}
\medskip

\noindent and, using the superadditivity of $\natural\,$, we get
$$
\aligned
f^\natural &\ge (1-\delta^2)\bigg(\sum^n_{m=1}
u^mww^*u^{-m}\bigg)^\natural\ge\\
&\ge(1-\delta^2)\bigg(\sum^n_{m=1}\big(
(u^mww^*u^{-m})^\natural\big)^k\bigg)^{\frac1k}=\\
&=(1-\delta^2)\Big( n\big( (w^*w)^\natural\big)^k\Big)^{\frac1k}
=(1-\delta^2)n^{\frac1k} p^\natural\, .
\endaligned
$$
But $p^\natural=e_j^\natural$ for all $1\le j\le n\,$, so
\medskip

\centerline{$n p^\natural=\sum^n_{j=1}e_j^\natural=
\left(\sum^n_{j=1}e_j\right)^\natural =1_M$}
\medskip

\noindent and we conclude that $f^\natural\ge(1-\delta^2)
n^{\frac1k-1}1_M\,$.

\hfill$\square$
\bigskip

Now we are ready to prove the following

\proclaim{Theorem 3 {\rm $($on the abelian strict closure in
continuous semi-finite $AW^*$-algebras$)$}} $\!$
Let $M$ be a continuous semi-finite \AW such that, for some
finite projection $e_o\in M$ of central support $1_M$ and
some $k\ge 1\,$, the centre valued quasitrace of $e_oMe_o$
is $k$-superadditive. Let further $A$ denote the norm-closed
linear span of all finite projections of $M\,$, and $C$
a masa of $A\,$. Then the strict closure of $C$ in $M=M(A)$
is $C\,$.
\endproclaim

{\bf Proof.} Let us assume that the strict closure
${\overline{C}}^\beta\subset C'\cap M$ of $C$ contains some
$0\le x\notin C\,$.
\smallskip

(a) First we prove that then ${\overline{C}}^\beta$ contains some
projection  $e\notin C\,$.

For let $e_\delta$ denote, for every $\delta>0\,$, the support
of  $(x-\delta1_M)_+$ in the \sAW $C'\cap M$ of $M\,$. Then
\medskip

\centerline{$xe_\delta\ge\delta e_\delta,\quad x(1_M-e_\delta)
\le\delta(1_M-e_\delta)\, .$}
\medskip

\noindent In particular, there exists $0\le y\in C'\cap M$
with $yx=e_\delta\,$. Moreover, $e_\delta\in{\overline{C}}^\beta\,$.
Indeed, by Lemma 1 (i) there is a net $(b_\iota)_\iota$ in $C$ with
\medskip

\centerline{$0\le b_\iota\le x\hbox{ for all } \iota,\quad
b_\iota\to x\hbox{ strictly.}$}
\medskip

\noindent Then $0\le yb_\iota\in C'\cap A=C$ for all $\iota$ and
\medskip

\centerline{$\|(e_\delta-yb_\iota)a\|=\|y(x-b_\iota)a\|\le
\|y\|\cdot\|(x-b_\iota)a\|\to 0\, .$}
\medskip

\noindent for every $a\in A\,$.
\smallskip

(b) Next we prove the existence of an infinite sequence of
mutually orthogonal projections $0\ne e_1,e_2,\,\dots\,\in C\,$,
all equivalent in $M$ to $e_oq_o$ for some projection $q_o$
in the centre $Z$ of $M\,$, such that $\bigvee_{n\ge 1} e_n
\in{\overline{C}}^\beta\,$.

Let $e$ be a projection as in (a). Then $e$ is not finite, so
there exists a projection $q\in Z$ such that $eq$ is properly
infinite. But then, by the comparison theorem, there exists a
projection $0\ne q_o\in Z$ such that $e_oq_o\prec eq\,$.
Since the central support of $e_o$ is $1_M\,$, we have $q_o\le q\,$.

Now, according to (ii) in Theorem 2 (on labeling Murray-von Neumann
equivalence classes), there exists a family $(e_\iota)_{\iota\in I}$
of mutually orthogonal projections in $C\,$, all equivalent in
$M$ to $e_oq_o\ne 0\,$, such that $\bigvee_{\iota\in I}e_\iota=
eq_o\, .\, I$ must be infinite, so it contains an infinite
sequence $\iota_1,\iota_2,\,\dots\;$. Put
\medskip

\centerline{$e_n=e_{\iota_n},\quad n\ge 1\, .$}
\medskip

\noindent Then $\bigvee_{n\ge 1}e_n$ belongs to
${\overline{C}}^\beta\,$. Indeed, since $\bigvee_{n\ge 1}e_n
\in C'\cap M\,$, if $(b_\kappa )_\kappa$ is a net in $C$ which
converges strictly to $e\,$, then the net $\big( b_\kappa
\bigvee_{n\ge 1}e_n\big)_\kappa$ is contained in $C$ and
converges clearly to $\, e\bigvee_{n\ge 1}e_n=\bigvee_{n\ge 1}e_n$
in the strict topology of $M\,$.
\smallskip

(c) Finally we prove that the statement in (b) leads to a
contradiction.

Let us denote by $\natural$ the map $\bigcup_{n\ge 1}e_nMe_n
\to Zq_o$ such that, for every $n\ge 1\,$, $e_nMe_n\ni x
\longmapsto x^\natural e_n$ is the centre valued quasitrace of
$e_nMe_n\,$. It is easy to see that $\natural$ takes the same
value in two projections from $\,\bigcup_{n\ge1}e_nMe_n$ if
and only if they are equivalent in $M\,$.

Let $n\ge 1$ be arbitrary and let $j_n=\big[ n^{\frac{k+1}k}
\big]\ge 1$ denote the integer part of $n^{\frac{k+1}k}\,$.
According to the corollary of Theorem 2 (on labeling
Murray-von Neumann equivalence classes), there exist projections
\smallskip

\centerline{$\displaystyle e_{n,1},\dots,e_{n,j_n}\in C,\quad
\sum^{j_n}_{j=1}e_{n,j}=e_n\, ,$}
\smallskip

\noindent such that
\smallskip

\centerline{$\displaystyle e_{n,j}^\natural =\frac1{j_n}q_o
\hbox{ for all } 1\le j\le j_n\, .$}
\medskip

\noindent Since $e_n\sim e_oq_o\,$, the centre valued
quasitrace of $e_nMe_n$ is $k$-superadditive and Lemma 4 yields
the existence of a projection $p_n\in e_nMe_n$ with
$p_n^\natural =\frac1{j_n}q_o$ such that, for every projection 
$g\in\{e_{n,1}\, ,\,\dots\, ,\, e_{n,j_n}\}'\cap e_nMe_n\,$,
\medskip

\centerline{$\displaystyle g^\natural\ge\big( 1-\|(e_n-g)p_n\|^2
\big)\frac1{j_n^{\frac{k-1}k}}\, q_o\ge\big( 1-\| (e_n-g)
p_n\|^2\big)\,\frac1n\, q_o\, .$}
\medskip

Now put $p=\bigvee_{n\ge 1}p_n\,$. Since $p_n^\natural =
\frac1{j_n}q_o$ and $\sum_{n\ge 1}\frac1{j_n}<+\infty\,$,
using Lemma 2 it is easy to verify that $p$ is equivalent to
a subprojection of the sum of finitely many $e_n$'s.
In particular, $p$ is finite, that is $p\in A\,$. Therefore,
$\bigvee_{n\ge 1}e_n$ being in ${\overline{C}}^\beta\,$,
Lemma 1 (iii) yields the existence of a projection
$\bigvee_{n\ge 1}e_n\ge f\in C$ with
\medskip

\centerline{$\displaystyle \bigg\|\bigg(\bigvee_{n\ge 1}e_n-
f\bigg)p\bigg\|\le\frac1{\sqrt 2}\, .$}
\smallskip

\noindent But then, for every $n\ge 1\, ,\, fe_n$ is a
projection in $C\cap e_nMe_n\subset\{e_{n,1}\, ,\,\dots\, ,\,
e_{n,j_n}\}'\cap e_nMe_n$ and the aboves yield
\medskip

\centerline{$\displaystyle (fe_n)^\natural\ge\big( 1-
\|(e_n-fe_n) p_n\|^2\big)\frac 1n\, q_o\ge\frac1{2n}\, q_o\, .$}
\medskip

\noindent Since $\sum_{n\ge 1}\frac1{2n} = +\infty\,$,
using again Lemma 2, it is easily seen that $f=\bigvee_{n\ge 1}
(fe_n)$ is equivalent to $\bigvee_{n\ge 1}e_n\,$. 
In particular, $f$ is properly infinite, in contradiction with
$f\in C\subset A\,$.

\hfill$\square$
\bigskip

{\bf 3 Weyl-von Neumann-Berg-Sikonia type theorems}
\bigskip

We recall that any Rickart \Cag $M$ is $\,\sigma$-normal, what
means that, for every increasing sequence $\big( e_k\big)_{k\geq 1}$
of projections in $M\,$, the least upper bound of
$\big( e_k\big)_{k\geq 1}\!$ in the projection lattice of $M$ is
actually its least upper bound in the ordered space $M_h$ of all
self-adjoint elements of $M$ (see [A-Go 2] or [Sa]). Therefore we
shall speak in the sequel simply about the least upper bound of
increasing sequences of projections in $M\,$.

Let us first prove a lemma about the sequential approximability
of a projection in a Rickart \Cag from a two-sided ideal :

\proclaim{Lemma 5} Let $M$ be a unital Rickart $C^*$-algebra,
$\Cal J$ a two-sided ideal of $M\,$, and $f\in M$ a projection.
Then the following statements are equivalent $:$
\smallskip

\itemitem{{\rm (a)}} $\,$there exists a sequence
$\big( b_k\big)_{k\ge 1}\!$ of positive elements in $\Cal J$ such
that $\, b_k\leq f$ for all $k\geq 1$ and every projection $e\in M$
with $\, b_k\leq e\, ,\, k\geq 1\,$, satisfies $f\leq e\,$;
\smallskip

\itemitem{{\rm (b)}} $\,$there exists an increasing sequence
$\big( f_k\big)_{k\geq 1}$ of projections in $\Cal J\,$, whose
least upper bound in $M$ is $f\,$.
\endproclaim

{\bf Proof.} Let us assume that (a) holds and put
\medskip

\noindent\hskip3.535cm $\displaystyle f_{k,l} =\hbox{ support of }
\Big( b_k-\frac1l\, 1_M\Big)_+\leq f\, ,\qquad k, l\ge 1\, ,$
\smallskip

\noindent\hskip1.29cm $\displaystyle f_n=
\bigvee_{1\le k,l\le n}f_{k,l}\,\hbox{ in the projection
lattice of }\, M\leq f\, ,\qquad n\ge 1\, .$
\smallskip

\noindent Since $\displaystyle b_k f_{k,l}\ge\frac1l f_{k,l}\,$,
and so $f_{k,l}$ can be factorized by $b_k\leq f\,$, we have
$f\geq f_{k,l}\in \Cal J$ for all $k$ and $l\,$.
Further, using the validity of the Parallelogramm Law in all Rickart
\Calg (see [Be], \S 13, Th. 1), we obtain also $f\geq f_n\in \Cal J\,
,\, n\geq 1\,$.

Now $(f_n)_{n\ge 1}$ is an increasing sequence, whose least
upper bound in the projection lattice of $M$ is $f\,$.
Indeed, if $e\in M$ is a projection which majorizes every $f_n\,$,
hence every $f_{k,l}\,$, then we have for all $k$ and $l$
\smallskip

\centerline{$\displaystyle b_k^{\frac12}\! (1_M-e)
b_k^{\frac12}\le b_k^{\frac12}\! (1_M-f_{k,l}) b_k^{\frac12}
\le \frac1l(1_M-f_{k,l})\, ,$}
\smallskip

\centerline{$\displaystyle \|(1_M-e)b_k^{\frac12}\|^2
\le\frac1l\, .$}
\smallskip

\noindent Thus
\smallskip

\centerline{$\displaystyle b_k=e\, b_ke\le e\hbox{ for all } k\ge 1$}
\medskip

\noindent and it follows that $f\leq e\,$.

Conversely, (b) obviously implies (a) with $b_k=f_k\,$.

\hfill$\square$
\bigskip

For unital Rickart \Calg we have the following
Weyl-von Neumann-Berg-Sikonia type result (cf. with [Z],
Theorem 3.1 and [Ak-Ped], \S 4) :
\medskip

\proclaim{Theorem 4}
Let $M$ be a unital Rickart $C^*$-algebra, and $\Cal J$ a
norm-closed two-sided ideal of $M\,$, which contains a sequence
of positive elements such that $1_M$ is the only projection in
$M$ majorizing the sequence. Then, for any normal $y\in M$ and
every $\ve >0\,$, there exist a masa $\, C$ of $\Cal J$ and an
element $x$ of the masa $\, C'\cap M$ of $M\,$, such that
\smallskip

\item\item{$1)$} $C$ contains an increasing sequence of projections,
whose least upper bound in

\noindent\hskip0.85cm $M$ is $1_M\,$,

\item\item{$2)$} $y-x\in \Cal J\,$ and $\,\|y-x\|\le \ve\,$.
\endproclaim
\medskip

\proclaim{Remark} We notice that in {\rm Theorem 4} $\, C'\cap M$
is the sequentially monotone closure of $\, C$ in the following
sense $:$ every $\, 0\leq a\in C'\cap M$ is the least upper bound
in $M_h$ of some increasing sequence of positive elements from
$\Cal J\,$.
\endproclaim

Indeed, if $\big( e_k\big)_{k\geq 1}$ is an increasing sequence of
projections in $\, C\,$, whose least upper bound in $M$ is $1_M\,$,
then $\big( a^{1/2}e_k a^{1/2}\big)_{k\geq 1}$ is an increasing
sequence of positive elements from $\Cal J\,$, whose least upper
bound in $A_h$ is $a^{1/2}1_M a^{1/2} =a$ (see [S-Z], 9.14, the
remark after Proposition 3).

\hfill$\square$
\medskip

For the proof of Theorem 4 we need the next result on quasi-central
approximate units, implicitely contained in [Z], Proposition 1.2 :
\medskip

\proclaim{Lemma 6} Let $M$ be a unital Rickart $C^*$-algebra,
$\Cal J$ an essential, norm-closed, two-sided ideal of $M\,$,
and $B\subset M$ a commutative $C^*$-subalgebra. Then the
upward directed set of all projections of $\Cal J$ contains
a subnet $(e_\iota)_{\iota\in I}$ which, besides being
automatically approximate unit for $\Cal J\,$, is quasi-central
for $B\,$, that is
\medskip

\centerline{$\displaystyle \lim_\iota\|\, e_\iota b-b\, e_\iota\|=
0\hbox{ for all } b\in B\, .$}
\endproclaim

{\bf Proof.} Passing to the Rickart \Csalg of $M$ generated
by $B$ and $1_M$ (see e.g. [S-Z], 9.11 (3)), we can assume
without loss of generality that $B$ is a Rickart \Csalg of $M$
containing $1_M\,$.

Let $\P$ denote the set of all finite sets $P$ of projections
from $B$ satisfying the equality $\sum_{p\in P} p=1_M$ and set 
\medskip

\centerline{$I=\{f\in A\, ; f\hbox{ projection }\} \times \P\, .$}
\medskip

\noindent We endow $I$ with a partial order by putting
$(f_1, P_1)\le (f_2, P_2)$ whenever $f_1\le f_2$ and the
\Cag $C^*(P_1)$ generated by $P_1$ is contained in $C^*(P_2)$
(that is the partition $P_2$ is a refinement of $P_1)\,$.
Clearly, in  this way $I$ becomes an upward directed ordered set.

Let $\iota=(f, P)\in I$ be arbitrary. For every $p\in P\,$, the
right support $\hbox{\bf r}(fp)$ of $fp$ is equivalent in $M$
to the left support $\hbox{\bf l}(fp)\le f\in \Cal J$ (see [A]
or [A-Go 1]), so it belongs to $\Cal J\,$. Thus
\smallskip

\centerline{$\displaystyle e_\iota=\sum_{p\in P}\hbox{\bf r}(fp)$}
\medskip

\noindent is a projection in $\Cal J\,$. Since every
$\hbox{\bf r}(fp)\le p$ belongs to the commutant $P'\,$, also
$e_\iota\in P'\,$. Furthermore,
\medskip

\centerline{$f\le e_\iota\, .$}
\medskip

\noindent Indeed, for every $q\in P,$
\smallskip

\centerline{$\displaystyle fq=fq\,\hbox{\bf r}(fq)=
\sum_{p\in P}fq\,\hbox{\bf r}(fp) =fq\, e_\iota\, ,$}
\smallskip

\noindent so
\smallskip

\centerline{$\displaystyle f=f\sum_{q\in P}q=
\sum_{q\in P}fqe_\iota=fe_\iota\le e_\iota\, .$}
\medskip

\noindent It is easily seen that
\smallskip

\centerline{$\iota_1\le \iota_2\Rightarrow e_{\iota_1}\le
e_{\iota_2}\, ,$}
\medskip

\noindent so $(e_\iota)_{\iota\in I}$ is a subnet of the
upward directed set of all projections of $\Cal J\,$.

Now, the upward directed set of all projections $f$ of $\Cal J$
is an increasing approximate unit for $\Cal J\,$. Indeed,
$\big\{ x\in \Cal J\, ;\,\lim\limits_f\|x(1_M-f)\|=0\big\}$ is
a norm-closed linear subspace of $\Cal J$ containing all
projections from $\Cal J\,$, hence it is equal to $\Cal J\,$.
Thus also the subnet $(e_\iota)_{\iota\in I}$ is an approximate
unit for $\Cal J\,$.

On the other hand, the norm-closed linear subspace
$\big\{ b\in B;\lim\limits_\iota\|e_\iota b-be_\iota\|=0\big\}$
contains every projection from $B\,$: for any
projection $p\in B$ and every $\iota=(f,P)$ with $p\in C^*(P)$
we have $e_\iota\in P'\cap A =C^*(P)'\cap \Cal J\,$, so
$e_\iota p-pe_\iota=0\,$. Consequently the above subspace of $B$
is actually equal to $B\,$.

\hfill$\square$
\bigskip

{\bf Proof of Theorem 4 .} Put $\displaystyle y_1=
\frac12\, (y+y^*)\, ,\, y_2=\frac1{2i}\, (y-y^*)$ and
\medskip

\centerline{$p_j(\lambda)=\hbox{ support of }(y_j-\lambda1_M)
\hbox{ in }M\, ,\quad\lambda\in \R\, .$}
\medskip

\noindent Let further $\{\lambda_1,\lambda_2,\,\dots\,\}$ be
the countable set of all rational numbers. Then
\smallskip

\centerline{$\displaystyle a=\sum^\infty_{k=1}3^{-(2k-1)}
(2p_1(\lambda_k)-1_{A^{**}})+\sum^\infty_{k=1}3^{-2k}
(2p_2(\lambda_k)-1_{A^{**}})+\frac121_{A^{**}} \in M\, ,$}
\medskip

\centerline{$0\le a\le1_M$}
\medskip

\noindent and it is easy to see that the \Csalg of $M$ generated
by $a$ and $1_M$ contains all projections $p_j(\lambda)\, ,\,
j=1,2\, ,\,\lambda\in\Bbb Q\,$, hence also $y=y_1+iy_2\,$.
Therefore there exists a continuous function $f : [0,+\infty)\to\C$
such that $y=f(a)\,$. Furthermore, by a well known continuity
property of the functional calculus (see e.g. [S-Z], 1.18 (5)),
there exists some $\delta>0$ such that
\medskip

\centerline{$0\le b\in M\, ,\, \|a-b\|\le\delta\;
\Longrightarrow\,\|f(a)-f(b)\|\le\ve\, .$}
\medskip

Now, by Lemma 5, there exists an increasing sequence
$\big( f_k\big)_{k\ge 1}$ of projections in $\Cal J\,$, whose
least upper bound in $M$ is $1_M\,$. Using Lemma 6, we can then
construct by induction a sequence $0=e_o\le e_1\le  e_2\le\dots$
of projections in $\Cal J$ such that
\medskip

\centerline{$f_k\le e_k\, ,\quad \|e_k a-ae_k\|\le
2^{-k-1}\delta\, .$}
\medskip

\noindent Since the elements $e_k$ and $(e_k-e_{k-1})a(e_k-e_{k-1})$
of $\Cal J$ are mutually commuting, there exists a masa $\, C$ of
$\Cal J$ containing all of them. Then $C$ contains the increasing
projection sequence $(e_k)_{k\ge 1}\,$, whose least upper bound in
$M$ is $1_M\,$\.

Let us denote
\medskip

\noindent\hskip1.257cm $b_o=a\, ,$
\smallskip

\noindent\hskip1.222cm $\displaystyle b_n=\sum^n_{k=1}
(e_k-e_{k-1})a(e_k-e_{k-1})+(1_M-e_n)a(1_M-e_n)\, ,
\quad  n\ge 1\, .$
\medskip

\noindent Then, for every $n\ge 1,$
\medskip

\noindent\hskip1.653cm $b_{n-1}-b_{n}=$
\smallskip

\noindent\hskip1.28cm $=(1_M-e_{n-1})\big( a-(1_M
-e_n)a(1_M-e_n) -e_nae_n\big) (1_M-e_{n-1})=$
\smallskip

\noindent\hskip1.28cm $=(1_M-e_{n-1})\cdot
[e_n\, ,\, e_na-ae_n]\cdot (1_M-e_{n-1})\, ,$
\medskip

\centerline{$\|b_{n-1}-b_n\|\le 2\|e_na-ae_n\|\le
2^{-n}\delta\, .$}
\medskip

\noindent It follows that $\sum\limits^\infty_{n=1}
\|b_{n-1}-b_n\|\le\delta\,$, so the sequence $(b_n)_{n\ge 1}$
is norm convergent to some $b\in M(A)^+$ and 
\medskip

\centerline{$\|a-b\|=\underset{n\to \infty}\to
\lim \|b_o-b_n\|\le\delta\, .$}
\smallskip

\noindent Put $x= f(b)\,$.

We claim that $b\in C'\cap M\,$, hence also $x\in C'\cap M\,$.
Since $C'\cap M$ is a masa of $M$ (see Lemma 1 (iv)), it is
enough to prove that $b$ is commuting with all elements
$a'\in C'\cap M\subset\big\{ e_k\, ,(e_k-e_{k-1}) a(e_k-e_{k-1})
\, ;\, k\ge 1\big\}'\cap M\,$. For we notice that, for every
$n\ge 1\,$,
\medskip

\centerline{$b_na'-a'b_n=(1_M-e_n)(aa'-a'a)(1_M-e_n)\, ,$}
\smallskip

\noindent hence
$$
\aligned
|b_na'-a'b_n|^2&\le (1_M-e_n)|aa'-a'a|^2 (1_M-e_n)\le \\
&\le\|aa'-a'a\|^2 (1_M-e_n)\, .
\endaligned
$$
Therefore
\medskip

\centerline{$|b_na'-a'b_n|^2\le\|aa'-a'a\|^2(1_M-e_k)
\, ,\quad n\ge k\ge 1$}
\medskip

\noindent and, passing to limit for $n\to\infty\,$,
we get for every $k\ge 1$
\medskip

\centerline{$|ba'-a'b|^2\le\|aa'-a'a\|^2 (1_M-e_k)\, ,$}
\medskip

\centerline{$\hbox{support of}\; |ba'-a'b|^2\;\hbox{in}\; M
\;\hbox{is} \le 1_M -e_k\, .$}
\medskip

\noindent Since the least upper bound of $(e_k)_{k\ge 1}$
in $M$ is $1_M\,$, it follows that $ba'-a'b=0\,$.

Finally, according to the choice of $\delta\, ,\,
\|a-b\|\le\delta$ implies that
\medskip

\centerline{$\|y-x\|=\|f(a)-f(b)\|\le \ve\, .$}
\medskip

\noindent On the other hand,
\medskip

\noindent\hskip1.76cm $\displaystyle a-b_n=\sum^n_{k=1}
(b_{k-1}-b_k)=$
\smallskip

\noindent\hskip2.818cm $\displaystyle =\sum^n_{k=1}
(1_M-e_{k-1})\cdot[e_k\, ,\, e_ka-ae_k]\cdot
(1_M-e_{k-1})\in \Cal J$
\medskip

\noindent implies by passing to the limit for $n\to\infty$
that $a-b\in \Cal J\, .$ Using the Weierstrass Approximation
Theorem, we infer that $y-x=f(a)-f(b)\in \Cal J\, .$

\hfill$\square$
\bigskip

We shall prove that in Theorem 4 the element $x$ can be found
under the form of an ``infinite linear combination'' of a
sequence of mutually orthogonal projections from $\Cal J\,$.
To this aim we need an appropriate understanding of the
summation of series in Rickart $C^*$-algebras.

We recall that every commutative Rickart \Cag $\, C$ is
sequentially monotone complete (see e.g. [S-Z], 9.16,
Proposition 1). Thus, if $\big( a_k\big)_{k\geq 1}$ is a
sequence in $C^+$ such that the partial sums
$\sum_{k=1}^n a_k\, , n\geq 1\,$, are bounded, then there
exists the least upper bound in $\, C_h$
\smallskip

\centerline{$\displaystyle \sum_{k=1}^\infty a_k =
\sup\bigg\{ \sum\limits_{k=1}^n a_k\, ;\, n\geq 1\bigg\}\in
C^+\, .$}
\smallskip

Let next $M$ be an arbitrary Rickart $C^*$-algebra,
$\big( a_k\big)_{k\geq 1}$ a bounded sequence in $M^+$ such
that the supports $\hbox{\bf s}(a_k)\,$, $k\geq 1\,$, are
mutually orthogonal, and $\big( e_k\big)_{k\geq 1}\!$ a
sequence of mutually orthogonal projections in $M\,$, for
which $\hbox{\bf s}(a_k)\leq e_k\, ,\, k\geq 1$ (we can take,
for example, $e_k=\hbox{\bf s}(a_k)$). Then $\big\{ a_k\, ;
\, k\geq 1\big\}\cup\big\{ e_k\, ;\, k\geq 1\big\}$ generates
a commutative Rickart $C^*$-subalgebra $\, C$ of $M\,$, so
there exists $a=\sum_{k=1}^\infty a_k\in C^+$.
Moreover, $a$ is the least upper bound of the partial sums
$\big\{\sum_{k=1}^n a_k\, ;\, n\geq 1\big\}$ even in $M_h\,$.
Indeed, by the $\sigma$-normality of the Rickart $C^*$-algebras,
$\bigvee_{k=1}^\infty e_k$ is the least upper bound in $M_h$
of the sequence $\big(\bigvee_{k=1}^n e_k\big)_{n\geq 1}$ and
it follows that
\smallskip

\noindent\hskip1.813cm $\displaystyle a=a^{1/2}\bigg(
\bigvee\limits_{k=1}^\infty e_k\bigg) a^{1/2}\,\hbox{ is the
least upper bound in }\, M_h\,\hbox{ of}$

\noindent\hskip1.813cm $\displaystyle \hbox{the increasing
sequence }\, a^{1/2}\bigg(\bigvee\limits_{k=1}^n e_k\bigg)
a^{1/2} =\sum\limits_{k=1}^n a_k\, ,\, n\geq 1$
\smallskip

\noindent (see [S-Z], 9.14, the remark after Proposition 3).
In particular, $a$ is the only element of $M_h$ satisfying
the conditions
\smallskip

\centerline{$\displaystyle a\, e_k=a_k\, ,\;\!
k\geq 1\, ,\qquad\hbox{\bf s}(a)\leq\bigvee\limits_{k=1}^\infty
e_k\, .$}
\smallskip

For sake of completeness we notice that, by the above
characterization, if $\big( e_k\big)_{k\geq 1}$ is a sequence
of mutually orthogonal projections in $M\,$, then
$\sum_{k=1}^\infty e_k =\bigvee_{k=1}^\infty e_k\,$.

Now let $\big( x_k\big)_{k\geq 1}$ be a bounded sequence in $M$
such that, denoting by $\hbox{\bf l}(x_k)$ the left support of
$x_k$ and by $\hbox{\bf r}(x_k)$ the right one, the projections
$\hbox{\bf l}(x_k)\vee\hbox{\bf r}(x_k)\,$, $k\geq 1\,$,
are mutually orthogonal. Then we can define
\smallskip

\centerline{$\displaystyle \sum_{k=1}^\infty x_k =
\Big(\sum_{k=1}^\infty (\Re x_k)_+ -\sum_{k=1}^\infty (\Re x_k)_-
\Big) +i\Big(\sum_{k=1}^\infty (\Im x_k)_+ -\sum_{k=1}^\infty
(\Im x_k)_-\Big)\, .$}
\smallskip

\noindent It is easy to see that, if $\big( e_k\big)_{k\geq 1}$
is any sequence of mutually orthogonal projections in $M$ such that
$\hbox{\bf l}(x_k)\vee\hbox{\bf r}(x_k)\leq e_k\, ,\, k\geq 1\,$,
then $\sum_{k=1}^\infty x_k$ is the only element $x\in M\,$, for
which

\noindent (*)\hskip2.075cm $\displaystyle x\, e_k=e_k x=x_k\, ,
\;\! k\geq 1\, ,\qquad\hbox{\bf l}(x)\vee\hbox{\bf r}(x)\leq
\bigvee\limits_{k=1}^\infty e_k\, .$
\smallskip

\noindent By the aboves, considering the direct product
$\, C^*$-algebra
\smallskip

\centerline{$\displaystyle \bigoplus_{k=1}^\infty\,
e_k M e_k =\bigg\{ \big( y_k)_{k\geq 1}\in\prod_{k=1}^\infty
e_k M e_k\, ;\, \sup\limits_{k\geq 1} \| y_k\| <+\infty\bigg\}\,$,}
\smallskip

\noindent the mapping

\centerline{$\displaystyle \bigoplus_{k=1}^\infty\,
e_k M e_k \ni\big( y_k)_{k\geq 1}\longmapsto\sum_{k=1}^\infty y_k
\in M$}
\smallskip

\noindent is well defined and it is an injective $*$-homomorphism.
Consequently
\smallskip

\noindent (**)\hskip3.98cm $\displaystyle \Big\|\,\sum_{k=1}^\infty
x_k\,\Big\| =\sup\limits_{k\geq 1} \| x_k\|\, .$

Finally, let $\big( e_k\big)_{k\geq 1}$ be a sequence of mutually
orthogonal projections in $M\,$, and $\displaystyle
\big( x_k)_{k\geq 1}\, ,\,\big( y_k)_{k\geq 1}\in
\bigoplus_{k=1}^\infty\, e_k M e_k\,$. Denoting by
$\overline{\hbox{lin}}\,\big\{ x_k-y_k\, ;\, k\geq 1\big\}$ the
norm-closed linear subspace of $M$ generated by $\big\{ x_k-y_k
\, ;\, k\geq 1\big\}\,$, we have
\smallskip

\noindent (***)\hskip1.035cm $\displaystyle \sum_{k=1}^\infty
x_k -\sum_{k=1}^\infty y_k\in\overline{\hbox{lin}}\,\big\{ x_k-y_k
\, ;\, k\geq 1\big\}\,$ if $\,\| x_k-y_k\|\longrightarrow
0\, .$
\smallskip

\noindent Indeed, according to (**), we have :
\smallskip

\centerline{$\displaystyle \Big\| \sum_{k=1}^\infty x_k -
\sum_{k=1}^\infty y_k -\!\!\!\!\!\!\!\!\!\!\!\!
\underbrace{\sum_{k=1}^n (x_k-y_k)}_{
\in\;\overline{\hbox{lin}}\,\big\{ {\textstyle x_k-y_k\, ;\,
k\geq 1}\big\}}\!\!\!\!\!\!\!\!\!\!\!\!\Big\| =
\sup\limits_{k\geq n+1}\| x_k-y_k\|\,
{\buildrel {n\to\infty}\over{\hbox to1.2cm{\rightarrowfill}}}
\,0\, .$}
\medskip

A slight modification of the proof of Theorem 4 yields the
following Weyl-von Neumann-Berg-Sikonia type result, which
is much closer to [Z], Theorem 3.1 than Theorem 4 :
\medskip

\proclaim{Theorem 5} 
Let $M$ be a unital Rickart $C^*$-algebra, and $\Cal J$ a
norm-closed two-sided ideal of  $M$, which contains a sequence
of positive elements such that $1_M$ is the only projection in
$M$ majorizing the sequence. Then, for any normal $y\in M$
and every $\ve >0\,$, there are
\smallskip

\noindent\phantom{xxx}-$\;\;$a sequence $(p_k)_{k\geq 1}$ of
mutually orthogonal projections in $\Cal J\,$,

\noindent\phantom{xxx}-$\;\;$a sequence $(\lambda_k)_{k\geq 1}$
in the spectrum $\,\sigma (y)\!$ of $\, y\,$,
\smallskip

\noindent such that
\smallskip

\item\item{$1)$} $\,$the least upper bound of $(p_n)_{n\geq 1}$
in $M$ is $1_M\,$,

\item\item{$2)$} $\displaystyle \, y-\sum_{k=1}^\infty \lambda_k
\, p_k\in \Cal J$ and $\,\displaystyle \,\Big\|\, y -
\sum_{k=1}^\infty \lambda_k\, p_k\Big\|\leq\ve\,$.
\endproclaim

{\bf Proof.} Repeating word for word the arguments from the first
paragraph of the proof of Theorem 4, we get $a\in M$ with
$0\le a\le 1_M$, a continuous function $f : [0,+\infty)\to\C$ and
$\delta >0\,$, such that $y=f(a)$ and
\medskip

\noindent $(\diamond )$\hskip2.168cm $0\le b\in M\, ,\, \|a-b\|
\le\delta\;\Longrightarrow\,\|f(a)-f(b)\|\le \ve\, .$
\medskip

Subtracting from $a$ an appropriate positive multiple of $1_M$
and modifying $f$ corrispondingly, if necessary, we can assume
that $0\in\sigma (a)\,$.

Choose a sequence $\delta /3 =\delta_1 >\delta_2 >\,\ldots\, >0$
which converges to $0\,$. According to the upper semicontinuity
of the spectrum, there exist
\medskip

\centerline{$\displaystyle \matrix
{} &{} &\eta_1\!\! &\!\! >\!\! &\!\!\eta_2\!\! &\!\! >\!\!
&\ldots &\!\! >\!\! & 0 \\
{} &{} &\wedge &{} &\wedge &{} &{} &{} &{} \\
\delta /3 &\!\! =\!\! &\!\!\delta_1\!\! &\!\! >\!\!
&\!\!\delta_2\!\! &\!\! >\!\! &\ldots &{} &{}
\endmatrix$}
\medskip

\noindent such that the spectrum of every $b\in M$ with $\| a-b\|
\leq\eta_k$ is contained in
\medskip

\centerline{$\displaystyle U_{\delta_k}\big(\sigma (a)\big)
=\big\{ \mu\in\C\, ;\, |\mu -\lambda (\mu )| <\delta_k
\,\hbox{ for some }\,\lambda (\mu )\in\sigma (a)\big\}\, .$}
\medskip

Arguing now again as in the proof of Theorem 4, we can
construct a sequence $0=e_o\le e_1\le  e_2\le\dots$ of projections
in $\Cal J\,$, whose least upper bound in $M$ is $1_M\,$, such that
\smallskip

\centerline{$\|e_k a-ae_k\|\le 2^{-k-1}\eta_{k+1}\,\hbox{ for all }\,
k\geq 1\, .$}
\medskip

\noindent Setting then
\medskip

\noindent\hskip1.257cm $b_o=a\, ,$
\smallskip

\noindent\hskip1.222cm $\displaystyle b_n=\sum_{k=1}^n
(e_k-e_{k-1})a(e_k-e_{k-1})+(1_M-e_n)a(1_M-e_n)\, ,
\quad  n\ge 1\, ,$
\medskip

\noindent we have
\medskip

\centerline{$b_{n-1}-b_n =(1_M-e_{n-1})\cdot
[e_n\, ,\, e_na-ae_n]\cdot (1_M-e_{n-1})\, ,\quad n\geq 1\, ,$}
\medskip

\noindent so $\| b_{n-1}-b_n\|\leq 2^{-n}\eta_{n+1}\leq 2^{-n}
\delta /3$ and $b_{n-1}-b_n\in \Cal J\,$. Therefore the sequence
$(b_n)_{n\ge 1}$ is norm convergent to some $b_\infty\in M^+$,
for which $\|a-b_\infty\|\le\delta /3$ and $a-b_\infty\in
\Cal J\,$.

We claim that
\smallskip

\centerline{$\displaystyle b_\infty=\sum_{k=1}^\infty (e_k-e_{k-1})a
(e_k-e_{k-1})\, .$}
\smallskip

\noindent Indeed, since
\medskip

\centerline{$\displaystyle b_n\, (e_k-e_{k-1}) =(e_k-e_{k-1})\, b_n
=(e_k-e_{k-1})a(e_k-e_{k-1})\, ,\qquad n\geq k\geq 1\, ,$}
\medskip

\noindent by passing to the limit for $n\to\infty$ we get
\medskip

\centerline{$\displaystyle b_\infty\, (e_k-e_{k-1}) =(e_k-e_{k-1})
\, b_\infty =(e_k-e_{k-1})a(e_k-e_{k-1})\, ,\qquad k\geq 1\, .$}
\medskip

\noindent Thus, taking into account that $
\bigvee_{k=1}^\infty (e_k-e_{k-1}) =\bigvee_{k=1}^\infty e_k =
1_M\,$, the description (*) yields the desired equality.

We notice that, for every $k\geq 1\,$,
\medskip

\noindent $(\diamond\diamond )$\hskip2.528cm
$\sigma\big((e_k-e_{k-1})a(e_k-e_{k-1})\big)\subset
U_{\delta_k}\big(\sigma (a)\big)\, .$
\medskip

\noindent Indeed, since the norm of
\smallskip

\centerline{$\phantom{xx} a-\Big( (e_k-e_{k-1})a(e_k-e_{k-1})
+\big( 1_{A^{**}}-(e_k-e_{k-1})\big) a\big( 1_{A^{**}}-
(e_k-e_{k-1})\big)\Big)$}

\noindent\hskip0.33cm $=\Big[ [e_k-e_{k-1},a]\, ,1_{A^{**}}-
(e_k-e_{k-1})\Big]$
\smallskip

\noindent is majorized by $2\;\!\big( \| e_ka-ae_k\| +
\| e_{k-1}a -ae_{k-1}\|\big)\leq 2\;\!\big( 2^{-k-2}\eta_{k+1}
+2^{-k-1}\eta_k\big) <\eta_k\,$, by the choice of $\eta_k$
we have
\medskip

\noindent\hskip0.431cm $\sigma\big((e_k-e_{k-1})a(e_k-e_{k-1})\big)$

\noindent\hskip0.058cm $\subset\sigma\Big( (e_k-e_{k-1})a
(e_k-e_{k-1}) +\big( 1_{A^{**}}-(e_k-e_{k-1})\big) a\big(
1_{A^{**}} -(e_k-e_{k-1})\big)\Big)\cup\{ 0\}$

\noindent\hskip0.058cm $\subset U_{\delta_k}\big(\sigma
(a)\big)\, .$
\medskip

For any $k\geq 1\,$, let $[r^{(k)}_1,r^{(k)}_2]$ denote
the smallest compact interval in $\R$ containing the spectrum
$\,\sigma\big((e_k-e_{k-1})a(e_k-e_{k-1})\big)\,$. Choose
\medskip

\centerline{$r^{(k)}_1 =\mu^{(k)}_1<\,\ldots\, <\mu^{(k)}_j<
\,\ldots\, <\mu^{(k)}_{j_k} =r^{(k)}_2$}
\smallskip

\noindent in $\,\sigma\big((e_k-e_{k-1})a(e_k-e_{k-1})\big)$
such that $|\mu^{(k)}_j -\mu^{(k)}_{j-1}|\leq\eta_k$ for al
$2\leq j\leq j_k\,$. Then there exist mutually orthogonal
projections $\big( p^{(k)}_j\big)_{1\leq j\leq j_k}$ in $\Cal J$
such that

\centerline{$\displaystyle \sum_{j=1}^{j_k}p^{(k)}_j=e_k-e_{k-1}\,$
and $\,\displaystyle \Big\| (e_k-e_{k-1})a(e_k-e_{k-1}) -
\sum_{j=1}^{j_k}\mu^{(k)}_j p^{(k)}_j\Big\|\leq\eta_k\, :$}
\smallskip

\noindent For example, we can set $p^{(k)}_j =
e^{(k)}_j -e^{(k)}_{j+1}\, ,\, 1\leq j\leq j_k\,$, where
\smallskip

\centerline{$e^{(k)}_j=\hbox{\bf s}\Big(\big( (e_k-e_{k-1})a
(e_k-e_{k-1}) -\mu^{(k)}_j\, (e_k-e_{k-1})\big)_+\Big)
\, ,\qquad 1\leq j\leq j_k$}
\smallskip

\noindent and $e^{(k)}_{j_k+1}=0$ (see e.g. [S-Z], 9.9, Proposition 1).
Using $(\diamond\diamond )$, we can find for every $\mu^{(k)}_j$ some
$\lambda^{(k)}_j\in\sigma (a)$ with $|\lambda^{(k)}_j -\mu^{(k)}_j|
<\delta_k$ and then
\smallskip

\centerline{$\displaystyle \Big\| (e_k-e_{k-1})a(e_k-e_{k-1}) -
\sum_{j=1}^{j_k}\lambda^{(k)}_j p^{(k)}_j\Big\|\leq\eta_k
+\delta_k <2\;\!\delta_k\leq 2\;\!\delta /3\, .$}

Now $\bigcup_{k=1}^\infty\big\{ p^{(k)}_j\, ;\, 1\leq j
\leq j_k\big\}$ consists of mutually orthogonal projections in
$M\,$, whose least upper bound in $M$ is $1_M\,$, while
$\bigcup_{k=1}^\infty\big\{ \lambda^{(k)}_j\, ;\, 1\leq j
\leq j_k\big\}\subset\sigma (a)\,$.
Set $\displaystyle b=\sum_{k=1}^\infty\sum_{j=1}^{j_k}\lambda^{(k)}_j
p^{(k)}_j\in M^+\,$. Then (**) yields

\noindent\hskip1.294cm $\displaystyle \| b_\infty -b\| =\Big\|
\sum_{k=1}^\infty (e_k-e_{k-1})a (e_k-e_{k-1}) -
\sum_{k=1}^\infty\sum_{j=1}^{j_k}\lambda^{(k)}_j p^{(k)}_j\Big\|$

\noindent\hskip2.774cm $\displaystyle =\,\sup\limits_{k\geq 1} \Big\|
(e_k-e_{k-1})a (e_k-e_{k-1}) -\sum_{j=1}^{j_k}\lambda^{(k)}_j
p^{(k)}_j\Big\|\leq 2\;\!\delta /3\, ,$
\medskip

\noindent so $\| a-b\|\leq \| a-b_\infty\| +\| b_\infty -b\|\leq
\delta /3 +2\;\!\delta /3 =\delta\,$.
On the other hand, since
\smallskip

\centerline{$\displaystyle \Big\| \underbrace{(e_k-e_{k-1})a
(e_k-e_{k-1})-\sum_{j=1}^{j_k}\lambda^{(k)}_jp^{(k)}_j}_{\in\,\Cal J}
\Big\| <2\;\!\delta_k\longrightarrow 0\, ,$}
\smallskip

\noindent (***) implies that $b_\infty -b\in\Cal J\,$, hence
$a-b =(a-b_\infty )+(b_\infty -b)\in \Cal J\,$.

Using the characterization (*), it is easy to deduce that
\smallskip

\centerline{$\displaystyle f(b)=\sum_{k=1}^\infty\sum_{j=1}^{j_k}
f(\lambda^{(k)}_j)\, p^{(k)}_j\, ,$}

\noindent where, by the Spectral Mapping Theorem,
$\bigcup_{k=1}^\infty\big\{ f(\lambda^{(k)}_j)\, ;\, 1\leq j\leq
j_k\big\}$ is contained in $f\big(\sigma (a)\big) =\sigma
\big( f(a)\big) =\sigma (y)\,$.
On the other hand, $(\diamond )$ yields the norm estimation
$\| y-f(b)\| =\| f(a)-f(b)\|\leq\ve\,$. Finally, using
$a-b\in \Cal J$ and the Weierstrass Approximation Theorem,
we infer also that $y-f(b) =f(a)-f(b)\in \Cal J\,$.

\hfill$\square$
\bigskip

If in the above theorem we are not requiring the norm estimation
in 2), then the coefficients $\lambda_k$ can be chosen even in the
essential spectrum of $y$ modulo $\Cal J\,$:
\medskip

\proclaim{Theorem 6} 
Let $M$ be a unital Rickart $C^*$-algebra, and $\Cal J$ a
norm-closed two-sided ideal of  $M$, which contains a sequence
of positive elements such that $1_M$ is the only projection in
$M$ majorizing the sequence. For any normal $y\in M$ there are
\smallskip

\noindent\phantom{xxx}-$\;\;$a sequence $(p_k)_{k\geq 1}$ of
mutually orthogonal projections in $\Cal J\,$,

\noindent\phantom{xxx}-$\;\;$a sequence $(\lambda_k)_{k\geq 1}$
in the spectrum $\sigma_{\Cal J}(y)$ of the canonical image of $y$
in the

\noindent\phantom{xxxxx}quotient $C^*$-algebra $M/\Cal J$
\smallskip

\noindent such that
\smallskip

\item\item{$1)$} $\,$the least upper bound of $(p_n)_{n\geq 1}$
in $M$ is $1_M\,$,

\item\item{$2)$} $\displaystyle \, y-\sum_{k=1}^\infty \lambda_k
\, p_k\in \Cal J\,$.
\endproclaim
\medskip

For the proof we need the next lifting result, which is
essentially [Z], Proposition 2.1 :
\medskip

\proclaim{Lemma 7} Let $M$ be a unital Rickart $C^*$-algebra,
and $\Cal J$ a norm-closed two-sided ideal of $M\,$.
For any self-adjoint $a\in M$ there exists a self-adjoint
$b\in M$ such that $\sigma (b) =\sigma_{\Cal J}(b)$ and
$a-b\in\Cal J\,$.
\endproclaim

{\bf Proof.} A moment's reflection shows that the proof of [Z],
Proposition 2.1 works for $M$ unital Rickart $C^*$-algebra
instead of $W^*$-algebra.

\hfill$\square$
\bigskip

{\bf Proof of Theorem 6.} Repeating again the arguments from the
first paragraph of the proof of Theorem 4, we get some $a\in M$
with $0\le a\le 1_M$ and a continuous function $f : [0,+\infty)
\to\C$ such that $y=f(a)\,$. Now, according to Lemma 7, there
exists a self-adjoint $b\in M$ such that $\sigma (b) =
\sigma_{\Cal J}(b)$ and $a-b\in\Cal J\,$. In particular,
$\sigma (b) =\sigma_{\Cal J}(a)\subset [0,1]\,$, and so
$0\le b\le 1_M\,$.

Let $x$ denote the normal element $f(b)\,$. Using the Weierstrass
Approximation Theorem, we infer that $y-x\in\Cal J\,$, hence, by
the Spectral Mapping Theorem, we have $\sigma (x)=
f\big(\sigma (b)\big) =f\big(\sigma_{\Cal J}(a)\big) =
\sigma_{\Cal J}(y)\,$. Now Theorem 5 yields the existence of
\smallskip

\noindent\phantom{xxx}-$\;\;$a sequence $(p_k)_{k\geq 1}$ of
mutually orthogonal projections in $\Cal J\,$,

\noindent\phantom{xxx}-$\;\;$a sequence $(\lambda_k)_{k\geq 1}$
in $\,\sigma (x) =\sigma_{\Cal J}(y)\,$,
\smallskip

\noindent such that the least upper bound of $(p_n)_{n\geq 1}$
in $M$ is $1_M$ and $\, x-\sum_{k=1}^\infty\lambda_k\, p_k\in
\Cal J\,$. Then $\, y-\sum_{k=1}^\infty\lambda_k\, p_k =
(y-x) +\big( x-\sum_{k=1}^\infty\lambda_k\, p_k\big)\in
\Cal J\,$.

\hfill$\square$
\bigskip

Let us say that a \Cag $A$ is $\sigma$-subunital if there exists
a sequence $(b_n)_{n\ge 1}$ in $A^+\,$, whose least upper bound
in $M(A)_h$ is $1_{A^{**}}\,$. Clearly, if $A$ is $\sigma$-unital
then it is $\sigma$-subunital. For commutative $A$ the two notions
coincide. However, if $M$ is a countably decomposable type
II${}_\infty$-factor and $A$ is the norm-closed linear span of
all finite projections of $M\,$, then $A$ is not $\sigma$-unital
(see [Ak-Ped], Prop. 4.5), but it is easily seen that it is
$\sigma$-subunital.

We remark that the sequence $(b_n)_{n\ge 1}$ in the definition
of the $\sigma$-subunitalness can be considered a kind of
``approximate unit with respect to the order structure''.
Indeed, according to [S-Z], 9.14, the remark after Proposition 3,
if the least upper bound of $(b_n)_{n\ge 1}$ in $M(A)_h$ is
$1_{A^{**}}$ and $x\in M(A)\,$, then the least upper bound of the
sequence $\big( x^*b_n x\big)_{n\ge 1}$ in $M(A)_h$ is $x^*x\,$.
\smallskip

By Theorems 5 and 6 we have :
\medskip

\proclaim{Corollary}
Let $A$ be a $\,\sigma$-subunital $C^*$-algebra, whose multiplier
algebra $M(A)$ is a Rickart $C^*$-algebra. For any normal
$y\in M(A)$ and any $\ve>0$ there exist
\smallskip

\noindent\phantom{xxx}-$\;\;$a sequence $(p_k)_{k\geq 1}$ of
mutually orthogonal projections in $A\,$,

\noindent\phantom{xxx}-$\;\;$a sequence $(\lambda_k)_{k\geq 1}$
in the spectrum $\,\sigma (y)\!$ of $\, y\,$,
\smallskip

\noindent such that
\medskip

\item\item{$1)$} $\,$the least upper bound of $(p_n)_{n\geq 1}$
in $M(A)_h$ is $1_{A^{**}}\,$,

\item\item{$2)$} $\displaystyle \, y-\sum_{k=1}^\infty \lambda_k
\, p_k\in A\,$ and $\,\displaystyle \,\Big\|\, y -
\sum_{k=1}^\infty \lambda_k\, p_k\Big\|\leq\ve\,$.

\noindent Moreover, if we don't require the second inequality in
$2)$, then the sequence $(\lambda_k)_{k\geq 1}$ can be chosen
even in the spectrum of the canonical image of $y$ in the corona
algebra $C(A)=M(A)/A\,$.
\endproclaim
\hfill$\square$
\bigskip

In particular, the above corollary can be applied to $A=K(H)\,$,
where $H$ is a separable complex Hilbert space, in which case the
series $\sum_{k=1}^\infty \lambda_k\, p_k$ converges even with
respect to the strict topology of $M(A)=B(H)\,$. This is the
statement of the classical Weyl-von Neumann-Berg-Sikonia Theorem,
but convergence with respect to the strict topology is used also
in its subsequent extensions to $\sigma$-unital $C^*$-algebras
with real rank zero multiplier algebra (see e.g. [M], [Br-Ped],
[Zh], [H-Ro], [L1], [L2], [L3]).

On the other hand, in the early extension from [Z] of the
Weyl-von Neumann-Berg-Sikonia Theorem to the norm-closed linear
span $A$ of all finite projections of an arbitrary semifinite
$W^*$-factor $M\,$, which for $M$ of type II${}_\infty$ turns
out to be not $\sigma$-unital, the series $\sum_{k=1}^\infty
\lambda_k\, p_k$ is proved to converge only with respect to the
$s^*$-topology. The reason, why here a weaker topology than the
strict topology should be used, is given by Theorem 3:
if $M$ is a type II${}_\infty$ $W^*$-factor and we assume that
a sum $\sum_{k=1}^\infty\lambda_k\, p_k$ with $p_k\in A$ is
strictly convergent, then, according to Theorem 3, we must have
$\sum_{k=1}^\infty\lambda_k\, p_k\in A\,$.
\bigskip

\newpage

\centerline{\bf Appendix}\par

\bigskip

We give here, for the convenience of the reader, a treatment of
a set-theoretical result of $T.$ Iwamura (see [Ma], Appendix II)
and two applications to the theory of $AW^*$-algebras.
\medskip

\proclaim{Proposition} Let $I,\le$ be an upward directed
partially ordered uncountable set. Then, there exist a well
order $\preccurlyeq$ on $I$ and a family $(I_\iota)_{\iota\in I}$
of subsets of $I$ such that

\hskip0.5cm$-$\hskip0.3cm $I_\iota$ is upward directed for
every $\iota\in I\, ,$

\hskip0.5cm$-$\hskip0.3cm {\rm card} $I_\iota <$ {\rm card} $I\,
,\quad \iota\in I\,$,

\hskip0.5cm$-$\hskip0.3cm $I_{\iota_1} \subset I_{\iota_2}$
whenever $\iota_1\prec\iota_2\, ,$

\hskip0.5cm$-$\hskip0.3cm $\bigcup_{\iota\in I} I_\iota=I\, .$
\endproclaim

{\bf Proof.} By Zermelo's theorem there exists a well order
$\preccurlyeq$  on $I\, .$ We can choose it such that
\medskip

\noindent (*)\hskip2.34cm $\hbox{card } \{\iota'\in I\, ;\,\iota'
\prec\iota\} <\hbox{ card }I\hbox{ for every }\iota\in I\, .$
\medskip

\noindent Indeed, if there exists some $\iota\in I$ such that
\medskip

\centerline{$\card\{\iota'\in I\, ;\,\iota'\prec\iota\}=
\card I\, ,$}
\medskip

\noindent then there exists a smallest $\iota$ with respect
to $\preccurlyeq\,$, having the above property. Choose for
this $\iota$ a bijection
\smallskip

\centerline{$\Phi : I\to \{\iota'\in I\, ;\,\iota'\prec\iota\}$}
\smallskip

\noindent and replace $\preccurlyeq$ by the well order,
according to which $\iota_1$ less or equal to $\iota_2$ means
$\Phi (\iota_1)\preccurlyeq\Phi(\iota_2)\,$.

We notice that, $I$ being infinite, (*) implies that $I$ does
not contain a largest element with respect to $\preccurlyeq\,$.

Let us denote
\medskip

\centerline{$J_\iota=\{\iota'\in I; \ \iota'\prec i \},\quad
\iota\in I\, .$}
\smallskip

\noindent Then
\smallskip

\centerline{$\hbox{card } J_\iota<\hbox{card } I\, ,
\quad \iota\in I\, ,$}
\smallskip

\centerline{$J_{\iota_1}\subset J_{\iota_2}\hbox{ whenever }
\iota_1 \prec\iota_2$}
\smallskip

\centerline{$\bigcup_{\iota\in I} J_\iota = I\, .$}
\medskip

\noindent On the other hand, $I,\le$ being upward directed,
we can choose for each finite $F\subset I$ some $\iota(F)\in I$
such that
\medskip

\centerline{$\iota\le\iota (F)\hbox{ for all } \iota\in F\, .$}
\medskip

\noindent Denote for every $J\subset I$
\medskip

\centerline{$D_1(J)=J\cup\{\iota(F)\, ;\, F\subset
J\hbox{ finite }\}\, .$}
\medskip

\noindent We notice that
$$
\aligned
&D_1(J)\hbox{ is finite  for $J$ finite, }\\
&\hbox{card } D_1 (J)=\hbox{card $J$ for $J$ infinite }
\endaligned
$$
and
\smallskip

\centerline{$D_1(J_1)\subset D_1(J_2)\hbox{ whenever } J_1
\subset J_2\, .$}
\medskip

Now we define by recursion
$$
\aligned
D_{n+1}(J)&=D_1(D_n(J))\supset D_n(J),\ n\ge 1\hbox{ integer,}\\
D_\omega(J)&=\bigcup_{n\ge 1}D_n(J)\, .
\endaligned
$$
Then

\centerline{$D_\omega(J)$ is countable for $J$ finite,}
\smallskip

\centerline{card $D_\omega(J)$ = card $J$ for $J$ infinite}

\noindent
and
\smallskip

\centerline{$D_\omega(J_1)\subset D_\omega(J_2)\hbox{ whenever }
J_1\subset J_2\, .$}
\medskip

\noindent Moreover, $D_\omega(J),\le$ is upward directed for
every $J\subset I\, .$

Now, putting
\medskip

\centerline{$I_\iota=D_\omega(J_\iota),\quad \iota\in I\, ,$}
\medskip

\noindent it is easy to see that all conditions from the
statement are satisfied.

\hfill$\square$
\bigskip 

The first corollary extends Lemma 3 (compare with [Be], \S 33,
Exercise 1):

\medskip

\proclaim{Corollary 1} Let $M$ be an $AW^*$-algebra, $f\in M$
a finite projection, and $(e_\iota)_{\iota\in I}$ an upward
directed family of projections in $M$ such that
\medskip

\centerline{$e_\iota\prec f\hbox{ for all } \iota\in I\, .$}
\smallskip

\noindent Then
\smallskip

\centerline{$\displaystyle \bigvee_{\iota\in I}e_\iota\prec f\, .$}
\endproclaim

{\bf Proof.} The case of countable $I$ can be easily reduced
to Lemma 6. Indeed, choosing a cofinal sequence $\iota_1\le
\iota_2\le\dots$ in $I\,$, we have
\bigskip

\centerline{$\displaystyle \bigvee_{\iota\in I}e_\iota=
\bigvee_{n\ge 1}e_{\iota_n}=e_{\iota_1}\lor
\bigvee_{n\ge 1}(e_{\iota_{n+1}}-e_{\iota_n})$}
\medskip

\noindent and we can apply Lemma 3 to $f$ and the family
$e_{\iota_1}\, ,\, e_{\iota_2}-e_{\iota_1}\, ,\,
e_{\iota_3}-e_{\iota_2}\, ,\,\dots\;$.

For the proof in the general case let $f\in M$ be a finite
projection and let us assume the existence of some upward
directed family $(e_\iota)_{\iota\in I}$ of projections in
$M$ such that
\medskip
\centerline{$\displaystyle e_\iota\prec f\hbox{ for all }
\iota\in I\, ,\hbox{ but }\bigvee_{\iota\in I}e_\iota
\nprec f\, .$}
\medskip

\noindent Choose among all such families one with $I$ of the
smallest cardinality. By the first part of the proof $I$ is
then uncountable.

Let the well order $\preccurlyeq$ on $I$ and the family
$(I_\iota)_{\iota\in I}$ of subsets of $I$ be as in the
above proposition.

According to the minimality property of card $I\,$, we have
\medskip
\centerline{$\displaystyle p_\iota=\bigvee_{\iota'\in I_\iota}
e_{\iota'}\prec f,\quad \iota\in I\, .$}
\medskip

\noindent On the other hand,
\medskip
\centerline{$p_{\iota_1}\le p_{\iota_2}\quad \hbox{whenever }
\iota_1\prec\iota_2\, ,$}
\medskip
\centerline{$\displaystyle \bigvee_{\iota\in I}p_\iota
=\bigvee_{\iota\in I}e_\iota\, .$}
\smallskip

\noindent Consequently, denoting
\medskip
\centerline{$\displaystyle q_\iota =p_\iota -
\bigvee_{\iota'\prec\iota}p_{\iota'}\le p_\iota\, ,\quad
\iota\in I\, ,$}
\smallskip

\noindent the projections $(q_\iota)_{\iota\in I}$ are
mutually orthogonal and
\medskip
\centerline{$\displaystyle \sum_{\iota\in F}q_\iota\prec f
\hbox{ for any finite } F\subset I\, .$}
\smallskip

\noindent By Lemma 6 it follows that
\smallskip
\centerline{$\displaystyle \bigvee_{\iota\in I}q_\iota\prec f\, .$}
\smallskip

\noindent But
\smallskip
\centerline{$\displaystyle \bigvee_{\iota\in I}q_\iota
=\bigvee_{\iota\in I} p_\iota=\bigvee_{\iota\in I}e_\iota\, .$}
\medskip

\noindent Indeed, otherwise it would exist a smallest
$\iota\in I$ with respect to $\preccurlyeq$ such that
$$
p_\iota\nleq\bigvee_{\iota'\in I}q_{\iota'}.
\leqno(**)
$$
But then we would have
\medskip
\centerline{$\displaystyle \bigvee_{\iota''\prec\iota}
p_{\iota''}\le\bigvee_{\iota'\in I}q_{\iota'}\, ,$}
\medskip

\noindent which contradicts $(**)\,$.

\hfill$\square$
\bigskip

For $M$ an arbitrary \AW and $Z$ a commutative \AW we call
\medskip
\centerline{$\Phi : \{ e\in M\, ;\, e\hbox{ projection }\}
\to Z^+$}
\medskip

\noindent normal if, for every upward directed family
$(e_\iota)_\iota$ of projections in $M\,$, we have
\medskip
\centerline{$\displaystyle \Phi\Big(\bigvee_\iota e_\iota\Big)
=\sup\Phi(e_\iota)\, ,$}
\smallskip

\noindent where sup denotes the least upper bound in
$Z^+\, .$ Clearly,
\medskip

\centerline{$\Phi$ normal $\Rightarrow \Phi$ completely additive,}
\medskip

\noindent but, using the above proposition similarly as
in the proof of the Corollary 2, we get also the converse
implication (which should be known, but for which we have
no reference):
\medskip

\proclaim{Corollary 2} Let $M\, ,\, Z$ be $AW^*$-algebras,
$Z$ commutative, and $\Phi : \{ e\in M\, ;\, e$ projection $\}
\to Z^+\, .$ Then
\medskip

\centerline{$\Phi$ normal $\Leftrightarrow \Phi$ completely
additive.}
\endproclaim
\smallskip

In particular, the centre valued dimension function of a
finite \AW is normal (see [Be], \S 33, Exercise 4). Also,
if $M$ is a discrete \AW and $e\in M$ is an abelian projection
of central support $1_M$, then the map $\Phi_e$ considered in
the proof of Theorem 1 (on the abelian strict closure in
discrete $AW^*$-algebras) is normal on the projection lattice
of $M\,$.
\bigskip\bigskip


\vglue 1 true cm
\par
Dipartimento di Matematica
\par 
Universit\`a di Roma "Tor Vergata"
\par
Via della Ricerca Scientifica
\par
00133 Roma - Italia
\par
E-mail addresses: 
\par
dantoni@axp.mat.uniroma2.it , zsido@axp.mat.uniroma2.it 

\end